\documentclass{article}

\pdfoutput=1
\usepackage{PRIMEarxiv}

\usepackage[utf8]{inputenc} % allow utf-8 input
\usepackage[T1]{fontenc}    % use 8-bit T1 fonts
\usepackage{hyperref}       % hyperlinks
\usepackage{url}            % simple URL typesetting
\usepackage{booktabs}       % professional-quality tables
\usepackage{amsfonts}       % blackboard math symbols
\usepackage{nicefrac}       % compact symbols for 1/2, etc.
\usepackage{microtype}      % microtypography
\usepackage{lipsum}
\usepackage{fancyhdr}       % header
\usepackage{graphicx}       % graphics
\graphicspath{{media/}}     % organize your images and other figures under media/ folder

\usepackage{multirow}
\usepackage{lscape}
\usepackage{caption}
\usepackage{subcaption}
\usepackage{rotating}
\usepackage{ctable}
\usepackage[vlined,linesnumberedhidden,resetcount,ruled]{algorithm2e} %
\usepackage{colortbl}
\usepackage{pifont}
\newcommand{\xmark}{\ding{55}}
\usepackage{natbib}
\usepackage{amsmath,amssymb,url,graphicx,color,array, amsthm}
\newcommand{\cT}{\mathcal{T}}
\newcommand{\cH}{\mathcal{H}}

\newtheorem{lemma}{Lemma}[section]
 \newtheorem{proposition}{Proposition}[section]

\title{Robust Direct Aperture Optimization for Radiation Therapy Treatment Planning
\thanks{\textit{\underline{Citation}}:  Danielle A Ripsman, Thomas G Purdie, Timothy CY Chan, and Houra Mahmoudzadeh. Robust Direct Aperture Optimization for Radiation Therapy Treatment Planning. \textit{INFORMS Journal on Computing}, 
2022.
\url{https://doi.org/10.1287/ijoc.2022.1167}}
}

\author{
  Danielle A. Ripsman \\
  Department of Management Sciences \\
  University of Waterloo \\
  ON, Canada\\
  \texttt{daripsman@uwaterloo.ca} \\
   \And
  Thomas G. Purdie \\
  Princess Margaret Cancer Centre \\
  Toronto \\
  ON, Canada\\
  \texttt{tom.purdie@rmp.uhn.on.ca} \\
   \AND
   Timothy C. Y. Chan \\
   Department of Mechanical and Industrial Engineering \\ University of Toronto \\
   ON, Canada \\
   \texttt{tcychan@mie.utoronto.ca} \\
   \And
   Houra Mahmoudzadeh \\
   Department of Management Sciences \\
 University of Waterloo \\
  ON, Canada\\
   \texttt{houra.mahmoudzadeh@uwaterloo.ca} 
}

\begin{document}
\maketitle

\begin{abstract}
Intensity-modulated radiation therapy (IMRT) allows for the design of customized, highly-conformal treatments for cancer patients. Creating IMRT treatment plans, however, is a mathematically complex process, which is often tackled in multiple, simpler stages. This sequential approach typically separates radiation dose requirements from mechanical deliverability considerations, which may result in suboptimal treatment quality. For patient health to be considered paramount, holistic models must address these plan elements concurrently, eliminating quality loss between stages. This combined direct aperture optimization (DAO) approach is rarely paired with uncertainty mitigation techniques, such as robust optimization, due to the inherent complexity of both parts. This paper outlines a robust DAO (RDAO) model and discusses novel methodologies for efficiently integrating salient constraints. Because the highly-complex RDAO model is difficult to solve, an original candidate plan generation (CPG) heuristic is proposed. The CPG produces rapid, high-quality, feasible plans, which are immediately clinically viable, and can also be used to generate a feasible incumbent solution for warm starting the RDAO model. Computational results obtained using clinical patient datasets with motion uncertainty show the benefit of incorporating the CPG, both in terms of first incumbent solution and final output plan quality.
\end{abstract}

% keywords can be removed

%%%%%%%%%%%%%%%~~~~~  SECTION 1 ~~~~~%%%%%%%%%%%%%%%
\section{Introduction}

According to estimates by GLOBOCAN, 2018 brought 18.1 million new cancer cases, and 9.6 million cancer deaths, globally \citep{bray2018global}. Radiation therapy (RT) is recommended over the course of the treatment process for roughly half of all cancer patients \citep{baskar2012cancer}. The goal of RT is to provide a treatment that eliminates a cancerous region (target) while sparing adjacent, healthy tissue. Today, many RT treatments are delivered using intensity-modulated RT (IMRT) which conforms plans to the targets using advanced mathematical planning \citep{bortfeld2006imrt}. 
\begin{figure}[htbp]
           \centering
       	\includegraphics[width=465px,height=62px]{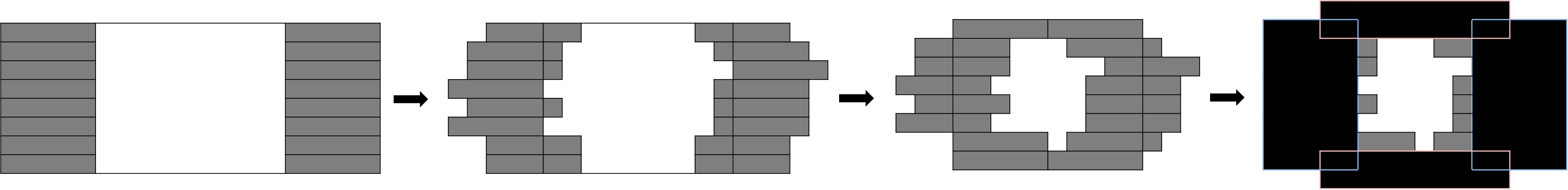} 
	 \caption{IMRT leaves move toward the centre of the beam to form an aperture. Peripheral parts of the leaves are then covered with insulating jaws.}
         \label{fig:MLC}
\end{figure} 

\begin{figure*}[htbp]
     \centering
     \hspace*{-2cm}
     \includegraphics[width=0.5\textwidth]{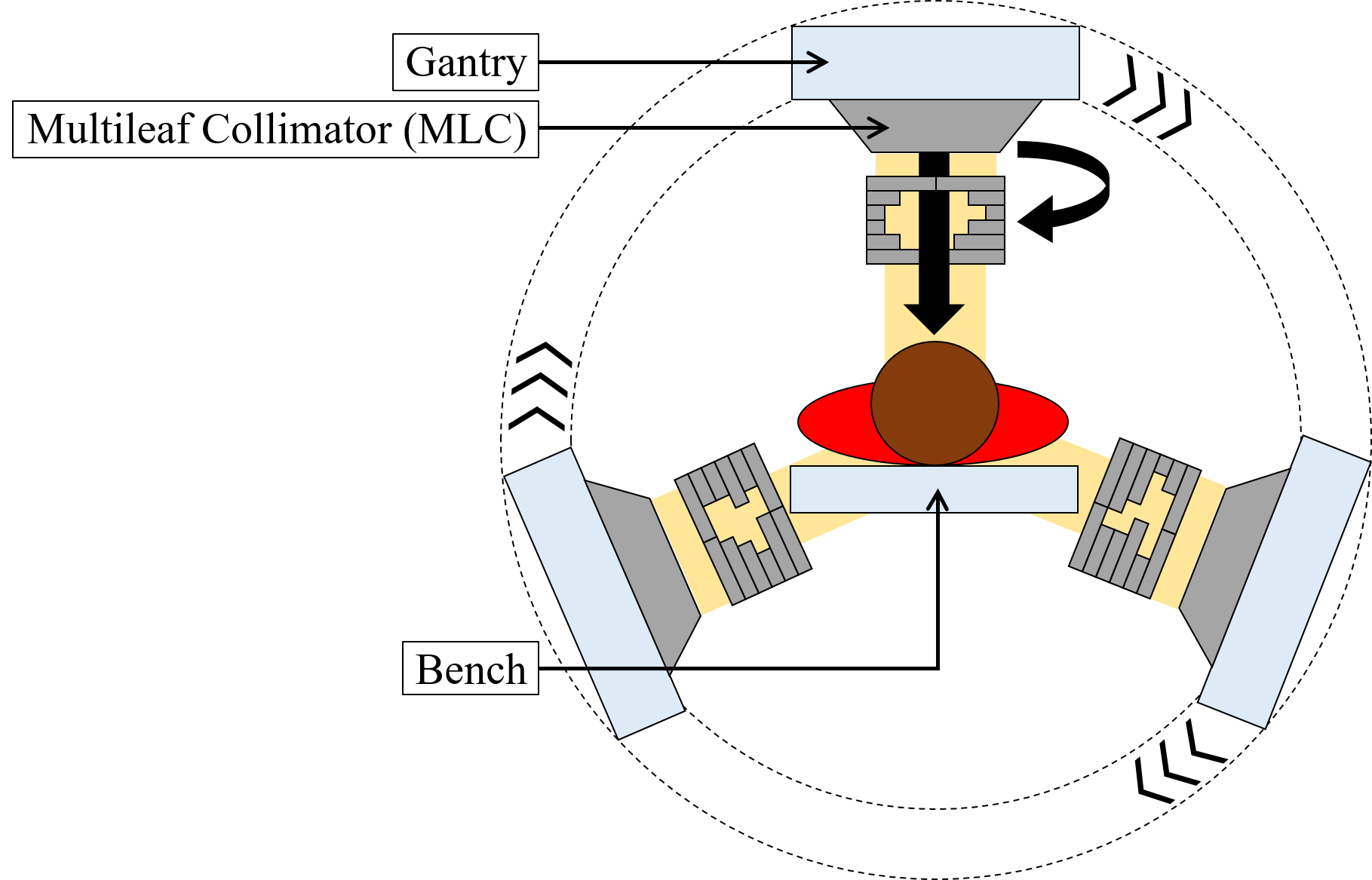} 
        \caption{The gantry rotates the MLC to various angles, delivering beams through a collection of apertures.}
         \label{fig:GenIMRT}
  \end{figure*}

IMRT is largely characterized by its modulating beam shapes, called \textit{apertures}. Apertures are formed when extending radiation-absorbing pieces of metal, called \textit{leaves} and \textit{jaws}, block out portions of an initially rectangular field during treatment, as shown in Figure~\ref{fig:MLC}.
In IMRT treatment delivery, the mouth of the delivery device, called a multileaf collimator (MLC), is rotated on an arm, called a gantry, to a specified angle where it delivers a set of shaped beams to the patient who is typically lying on a bench, as shown in Figure~\ref{fig:GenIMRT}. 
  
An IMRT plan consists of a set of predetermined apertures, each with an associated beam angle and fluence intensity. For simplicity, this work will focus on \textit{step-and-shoot} treatment plans, wherein the apertures are delivered statically and mechanical re-orientations or \textit{steps} between apertures occur with the beam turned off. The step-and-shoot planning techniques discussed may be generalized to commonly-used dynamic treatments, such as Volumetric Modulated Arc Therapy (VMAT), with the placement of additional fluence intensity and aperture shape requirements \citep{craft2014plan}.
 
\begin{figure}[htbp]
           \centering
       	\colorbox{white}{\includegraphics[width=.8\textwidth]{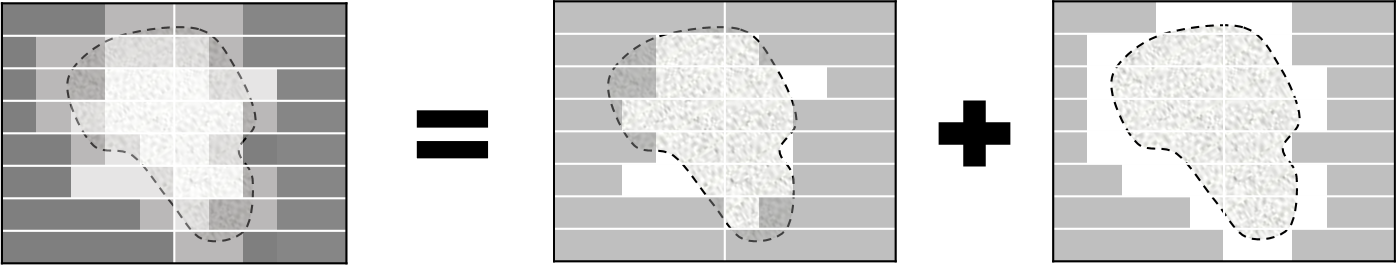}} 
	 \caption{Single-angle Step~1 fluence map broken down into two deliverable apertures in Step~2.}
         \label{fig:fmo}
\end{figure}

Even when plans are static, calculating the optimal combination of beam angles, apertures and intensity values simultaneously is an intractably large combinatorial version of an already NP-hard problem \citep{sultan2006optimization}. 
Although beam angle selection is an active area of research, empirical selection methods, such as choosing a prespecified number of equidistant angles, are often used in practice 
\citep{jiang2005examination,bertsimas2013hybrid}.
Once beam angles are chosen, the selection of apertures and their intensity values is often structured as a two-step problem:
\begin{enumerate}
\item[]\textbf{Step 1.} Solve a fluence map optimization (FMO) problem.
\item[]\textbf{Step 2.} Decompose the FMO into finite deliverable apertures and intensities.
\end{enumerate}
\noindent \textbf{Step~1} has received the bulk of the attention in the literature over the years, as the FMO model serves as a base for studying most clinically-significant objectives. While some of these objectives are inherently non-linear, many have been shown to have analogous linearizations, allowing for linear programming (LP) formulations of the FMO problem that can be rapidly solved to optimality \citep{romeijn2003novel, ehrgott2010mathematical, chan2014robust}. Regardless of formulation, during FMO, the beam is typically divided into a finite grid and each fluence in the grid is chosen independently, neglecting real-life delivery constraints. The left-hand side of Figure~\ref{fig:fmo} shows a hypothetical FMO output, which is optimal for target coverage but undeliverable for an MLC.

\noindent \textbf{Step~2}, depicted on the right-hand side of Figure~\ref{fig:fmo}, is a secondary decomposition model that follows the FMO problem \citep{bortfeld2006imrt,tacskin2010optimal,jing2015effective}. 
The decomposition often involves a heuristic leaf-sequencing algorithm, due to the mathematical difficulty associated with finding a clinically viable breakdown of an FMO solution \citep{carlsson2008combining}. Some studies segment the fluence maps into discrete fluence levels initially and then apply sequencing algorithms \citep{gladwish2007segmentation, xia1998multileaf, kamath2003leaf}. Others contribute exact algorithms for leaf sequencing when the goal is minimizing beam-on-time (BOT), i.e., the total duration of radiation delivery \citep{siochi1999minimizing, langer2001improved}. The leaf sequencing problem with a BOT minimization objective is structurally similar to network flow problems, a feature that has been exploited to create polynomially solvable LP decomposition models \citep{boland2004minimizing, ahuja2004network, tacskin2010optimal}. However, choosing more realistic objectives for this step of the problem, including minimizing the total number of apertures or total treatment time, leads to NP-hard problems \citep{baatar2005decomposition}.

When combined, the two steps of traditional step-and-shoot planning form a computationally challenging problem called direct aperture optimization (DAO). The DAO approach was initially introduced using a simulated annealing algorithm to generate the uniform apertures \citep{shepard2002direct}. This work was later augmented by a number of inexact methods for generating direct aperture plans \citep{li2003genetic, milette2008direct, broderick2009direct}. Efforts to find a globally optimal set of apertures have led to column generation approaches for both DAO 
and the dynamic VMAT version of the problem,
which can be solved 
(or approximately solved)
quite rapidly \citep{romeijn2005column, men2007exact, carlsson2008combining, salari2013column, mahnam2017simultaneous, dursun2019using}. 
While the column generation approaches in the literature are fast to converge, they do not allow for a hard limit on the number of apertures to be imposed \textit{a priori}, which may translate to unreasonably long total treatment times and negligibly small apertures, or a loss in optimality, when it comes to step-and-shoot planning. 
A framework for creating optimal treatment plans using a constrained number of apertures requires the explicit modeling of clinical and practical constraints in a single optimization model; a task 
that was recently examined in a dynamic setting \citep{akartunali2015unified, dursun2016mathematical}, but
has not previously been considered in static planning,
to the best of our knowledge.

In addition to delivery issues, theoretically optimal treatment plans may suffer from data uncertainty, which can adversely affect the quality of IMRT treatment in practice. Robust FMO (i.e., RFMO, the robust formulation of Step~1) is a popular method for immunizing problems against the parameter uncertainty seen in many forms of RT planning, from mechanical setup to patient and organ location uncertainty \citep{bertsimas2011theory,unkelbach2018robust}. Using robust optimization to handle RT uncertainty for patient setup and organ motion was first proposed by \cite{chu2005robust}. Similarly, uncertainties in dose calculation \citep{olafsson2006efficient}, machine parameter and delivery uncertainty, and machine parameter and organ location uncertainty \citep{bertsimas2010nonconvex,cromvik2010robustness}, have all been addressed using robust optimization. Intra-fractional breathing motion has also been captured by RFMO, where the uncertainty set is the proportion of time spent in each phase of a breathing cycle from inhale to exhale \citep{chan2006robust,bortfeld2008robust}. Continued studies by \cite{chan2014robust} and \cite{mahmoudzadeh2015robust, mahmoudzadeh2016constraint} have further solidified the robust model as a flexible and extensible model for addressing phase-based uncertainty in RT planning. 

Robust optimization in conjunction with DAO (RDAO) has been used to mitigate the impact of the tongue-and-groove effect in IMRT devices; a phenomenon that occurs when the interlocking MLC leaf edges are left exposed in an aperture, blocking radiation intended to reach the target \citep{salari2011accounting}. However, other practical phase-based sources of uncertainty, such as motion uncertainty, have not been considered in the presence of DAO. Some studies have covered the impact of motion uncertainty in a non-robust capacity using a commercial planning software \citep{zhang2006direct}. Alternatively, \citet{ahunbay2007investigation} account for motion using a gating system along with a commercial DAO planning system.
To the best of our knowledge, however, there is no unified approach for embedding more generalized forms of robust optimization, into a complete RDAO model. 

This paper presents a holistic mixed integer programming (MIP) framework for combining robust and DAO modeling methodologies. Due to the computational difficulty of solving this large-scale MIP model, novel symmetry-breaking constraints are introduced for angle allocation, as well as a more practical candidate plan generation (CPG) heuristic for rapidly finding viable plans. These CPG plans can be used as stand-alone treatments, or converted into a basic feasible solution for the RDAO model. It should be noted that the focus of this paper is providing a mathematical proof-of-concept for mechanical deliverability constraints, with little focus on advanced dose and clinical constraints commonly discussed in the literature.
More formally, the contributions of this paper include:
\begin{itemize} 
\item The development of a MIP framework for integrating robust and DAO constraints. 
\item An original formulation for five types of
mechanical deliverability constraints.
\item An exploration of symmetry-breaking techniques for aperture-angle assignments.
\item A heuristic for rapidly generating high-quality, feasible, robust DAO treatments.
\item An algorithm for converting feasible robust DAO treatments into initial incumbent solutions for warm starting commercial solvers.  
\end{itemize} 

The rest of this paper is organized as follows. Section~\ref{Sec:formulation} builds up the holistic model, starting with a basic RFMO model, then introducing and integrating compatible DAO constraints. Section~\ref{Sec:warmStart} proposes the CPG heuristic for creating deliverable plans, followed by an algorithm for converting these plans into warm starts for the RDAO model. Section~\ref{Sec:Results} provides the results of running the proposed RDAO model and CPG heuristic (both separately and together) on clinical breast cancer case studies with breathing motion uncertainty.
Finally, Section~\ref{Sec:Conclusion} contains concluding remarks, limitations and future research directions.

%%%%%%%%%%%%%%%~~~~~  SECTION 2 ~~~~~%%%%%%%%%%%%%%%

\section{Problem Formulation}
\label{Sec:formulation}
The RT treatment planning problem can be defined in terms of the relationship between the \textit{beam} of radiation and the patient's \textit{region(s) of interest}. 
Under this framework, the beam is broken down into a grid of $b \in \mathcal{B}$ rectangular units called beamlets. 
A region of interest is comprised of a set of $s \in \mathcal{S}$ structures and can be modeled as a finite grid of three-dimensional (3D) units called voxels, which are denoted by $v \in \mathcal{V}_s$. Typically, the structure set $\mathcal{S}$ is partitioned into targets $t \in \mathcal{T}$ and healthy organs $h \in \mathcal{H}$. 

The dosimetric influence that each beamlet $b$ has on each voxel $v$, per unit of radiation intensity, is captured by a \textit{dose influence matrix} $D_{v,b}$.
In standard RT treatment planning, the dose influence matrix is constant throughout the treatment.
One way to capture changes or uncertainty in organ shape, location or dose delivery during treatment is to discretize the changes into phases, allowing the dose influence matrix to take on a third \textit{phase} dimension.
Under breathing motion uncertainty, for example, each phase $i$ is associated with a snapshot of the healthy and target organs at a specified point in time.
Adjusted parameter $D_{v,b,i}$ captures the influence that a unit of intensity of each beamlet $b$ has on each voxel $v$ during phase $i$. Each phase $i \in \mathcal{I}$ also has a corresponding parameter $p_i$, which is the proportion of time spent in phase $i$, over the course of a complete treatment cycle. 

In this section, a complete robust direct aperture optimization model for IMRT planning is outlined. 
Section~\ref{Sec:baseMod} provides the base, robust fluence map optimization model, ignoring deliverability constraints.
Section~\ref{Sec:Direct} introduces five novel formulations for direct aperture optimization constraints. Finally, Section~\ref{sec:abridgedMod} contains a summary of the proposed robust direct aperture optimization model. Note that calligraphic capital letters will continue to denote sets throughout this paper, whereas bold letters indicate vectors or matrices.

\subsection{Basic Robust FMO Model}
\label{Sec:baseMod}
In a typical fluence map optimization (FMO) model, each beamlet $b$ may assume any intensity value, $\omega_b$, independent of its neighboring beamlets.  While there are many ways to model this problem (see~\cite{romeijn2008intensity} for a discussion about selecting clinical goals and objectives) this paper will use FMO to refer to the following simple LP problem. The objective concerns minimizing the expected dose to healthy organ voxels $v \in \mathcal{V}_\mathcal{H}$, while controlling for target voxel $v \in \mathcal{V}_\mathcal{T}$ overdose. Constraints ensure a prescribed dose is met, the simplest form of which enforces a uniform level of dose $L_v$ in each target voxel.
Finally, the model becomes a robust FMO (RFMO) when it addresses some worst-case parametric uncertainty e.g., the proportion of time spent in each phase $i\in \mathcal{I}$. Uncertain parameter $\tilde{p_i}$ may fall in uncertainty set $\mathcal{P}$, resulting in the following model. 
\begin{subequations}
\begin{align}
\tag{\textbf{RFMO}} \label{RobMod}\\
\min
& \sum_{s\in \{\cT,\cH\}}\frac{c_s}{|\mathcal{V}_s|}\sum_{v \in \mathcal{V}_s} \sum_{b \in \mathcal{B}} \sum_{i \in \mathcal{I}} p_i D_{v,b,i} \omega_{b} \label{RobObj} \\% \nonumber \\
\text{s.t.} & ~
 \sum_{b \in \mathcal{B}} \sum_{i \in \mathcal{I}} \tilde{p}_i D_{b,v,i} \omega_{b} \geq L_{v} & \forall v \in \mathcal{V}_T, \tilde{\textbf{p}} \in \mathcal{P}, \label{RobModInf}\\ 
& \omega_{b} \geq 0  & \forall b \in \mathcal{B}, \label{RobPos}
\end{align}
\end{subequations}
where $c_s$ is the objective weight for structure $s\in \mathcal{S}$. For a nominal FMO model, the robustness assumption is dropped, and fixed parameter $p_i$ dictates the fixed proportion of time spent in each phase $i\in \mathcal{I}$.
Static (single-phase) treatments can be captured by \eqref{RobMod} as well, by setting $|\mathcal{I}| = 1$.  

Note that the complexity of \eqref{RobMod} depends on the shape of the uncertainty set $\mathcal{P}$. If for instance, $\mathcal{P}$ is ellipsoidal, there exists a robust counterpart that is a second-order conic program \citep{ben1999robust}, whereas the following polyhedral set often seen in RT,
\begin{equation} 
\label{PDef} 
\mathcal{P} = \{\tilde{\textbf{p}} \in \mathbb{R}^{|\mathcal{I}|} | ~ p_i - \underline{p}_i \leq \tilde{p}_i \leq p_i + \bar{p}_i ~ \forall i \in \mathcal{I}; \sum_{i \in \mathcal{I}} \tilde{p}_i = 1; ~ 0 \leq \tilde{p}_i \leq 1 ~  \forall i \in \mathcal{I} \}, 
\end{equation}
has been shown by \cite{chan2006robust} to lead to the linear program in 
%Online
Appendix \ref{Sec:robMod}. 
%A.
The set \eqref{PDef} places bounds on each element in vector $\tilde{\textbf{p}}$ with a set of upper and lower deviations from the nominal proportions, denoted as $\bar{\textbf{p}}$ and \textbf{$\underline{\text{p}}$}, respectively, as well as enforcing feasible cycles with the $\sum_{i \in \mathcal{I}} \tilde{p}_i = 1$ term.

Also note that \eqref{RobMod} objective function, \eqref{RobObj}, does not contain uncertain parameters due to an RT practice called fractionation. Fractionation is the division of a treatment plan into smaller doses called ``fractions'', which are then delivered over multiple treatment sessions. This practice leads to any uncertainty in the objective function largely washing out over a complete treatment, while the constraints must hold at each individual fraction to ensure complete eradication of target cells \citep{bortfeld2008robust}.

\subsection{DAO Constraints}
\label{Sec:Direct}
In this section, we propose additional constraints to transform \eqref{RobMod} from a continuous FMO, to a linear, mixed integer model that produces distinct apertures, deliverable by IMRT collimators. Prior to adding constraints, the decision variable $\omega_b$ must be adjusted to handle direct aperture optimization (DAO) requirements. To this end, a new decision variable $w_{b,a}$ is introduced, which captures both beamlet intensity and the assignment of those intensities to one of the $a \in \mathcal{A}$ apertures. Substituting,
\begin{align}
\omega_b = \sum_{a \in \mathcal{A}} w_{b,a} ~~ \forall b \in \mathcal{B}, \label{sumW}  \end{align}
makes \eqref{RobMod} DAO-compatible, meaning it may be integrated with the proposed DAO constraints in the upcoming subsections, without any additional modification. The following subsections introduce five sets of constraints that ensure model intensity output is compatible with conventional IMRT multileaf collimators (MLCs). Section~\ref{Sec:unifCons} enforces aperture fluence uniformity, Section~\ref{Sec:Aperture} 
governs aperture-angle allocation and Section~\ref{Sec:Island} prevents MLC leaf discontinuities.
Finally, Sections~\ref{sec:VertCont} and \ref{sec:HorCont} provide optional constraints to rule out undesirable vertical and horizontal aperture breaks, respectively. 

\subsubsection{Uniformity Constraints.}
\label{Sec:unifCons}
Each active beamlet in a given aperture $a$ must either take on a fixed intensity value (to be chosen by the model), or have no intensity at all. Mathematically, this requirement may be enforced in \textbf{w} through the introduction of binary variables $x_{b,a}$ and continuous aperture intensity variable $f_a$. The $x_{b,a}$ variables indicate whether beamlet $b$ in aperture $a$ is open (i.e., no leaf is blocking the radiation) or shut (i.e., blocked by a leaf), as follows,
\[  w_{b,a}  = 
 \left\{
\begin{array}{ll}
      f_a & \text{if beamlet $b$ is active } (x_{b,a} =1),  \\
      0 & \text{if beamlet $b$ is inactive } (x_{b,a}=0). \\
\end{array} 
\right. \]
This property is achieved by adding the following constraints to the model:
\begin{subequations}
\label{Cons:uniform}
\begin{align}
&w_{b,a}  \leq Mx_{b,a} && \forall b \in \mathcal{B}, a \in \mathcal{A}, \label{blockedBeam}  \\ 
&w_{b,a}  \leq f_a + M(1 - x_{b,a}) && \forall b \in \mathcal{B}, a \in \mathcal{A}, \label{openBeam}  \\ 
&w_{b,a}  \geq f_a -M(1 - x_{b,a}) && \forall b \in \mathcal{B}, a \in \mathcal{A}, \label{openBeam2}\\ 
&f_a  \geq 0 && \forall a \in \mathcal{A}, \label{PositiveDose}\\
&x_{b,a} \in \{0,1\}  && \forall b \in \mathcal{B}, a \in \mathcal{A}, \label{binCon}
\end{align}
\end{subequations}
where $M$ is a very large number that may be selected based on the machine-specific dose rate and a clinically-determined upper bound on allowable treatment time. 

\subsubsection{Aperture Selection.}
\label{Sec:Aperture}
The set $\mathcal{B}$ includes all beamlets over all angles 
$\theta \in \Theta$. Let $\mathcal{B}_1, \dots, \mathcal{B}_{|\Theta|}$ be a partitioning of the set of all beamlets $\mathcal{B}$, where $\mathcal{B}_{\theta}$ is the set of beamlets at angle $\theta$. 
This implicit beamlet partitioning, defined by $b \in \{ 1,\dots,|\mathcal{B}|\} = \{1,\dots,|\mathcal{B}_1|,  |\mathcal{B}_1|+1,\dots,|\mathcal{B}_2|,\dots,|\mathcal{B}_{|\Theta|}| \}$, 
is not accounted for in uniformity constraints \eqref{Cons:uniform},
which may lead to single apertures that span multiple beam angles in a given treatment plan.
In order to generate deliverable plans, non-zero fluence intensities must be restricted to a single angle at each aperture. 

One way to enforce the angle-aperture allocation is through \textit{preallocation}, i.e., assigning a prespecified number of apertures to each angle, at a potential loss of global optimality. Alternatively, the decision to activate any single angle $\theta$ for a given aperture $a$, can be incorporated into a \textit{decision-based} model using binary decision variable $u_{a,\theta}$, where,
\[  u_{a,\theta}  = 
 \left\{
\begin{array}{ll}
      1 & \text{if angle $\theta$ is active in aperture $a$},   \\
      0 & \text{if angle $\theta$ is inactive in aperture $a$}. \\
\end{array} 
\right. \]
The assignment is then enforced through the following constraints,
\begin{subequations}
\label{Cons:fullApt}
\begin{align}
&\sum_{b \in \mathcal{B}_\theta} x_{b,a} \leq |\mathcal{B}_\theta|u_{a,\theta} && \forall a \in \mathcal{A}, \theta \in \Theta, \label{angleChoice}\\
&\sum_{\theta \in \Theta} u_{a,\theta} = 1 && \forall a \in \mathcal{A}. \label{oneAngle}
\end{align}
\end{subequations}
This decision-based allocation increases model flexibility, but it comes with a cost, both in terms of model size and symmetry in the solution space, as shown in Lemma \ref{lem:symGrowth}. Note that in the context of this section, preallocated models differ from decision-based models only in that binary decision variable
\textbf{u} in constraints~\eqref{Cons:fullApt} is replaced with a matrix of \textit{a priori}-determined
$\{0,1\}$ parameters.
\begin{lemma}
\label{lem:symGrowth}
Let $P_d$ and $P_p$ be the number of aperture permutations in decision-based and preallocated models, respectively.
 If $u_{a_1,\theta_1} = 1$ and $\exists u_{a\neq a_1,\theta\neq\theta_1} = 1$, then $P_d > P_p$.
\end{lemma}
\proof%{Proof.}
When aperture allocation is decision-based, there are $P_{d} = |\mathcal{A}|!$ ways to sequence the set among $|\Theta|$ angles. If $|\mathcal{A}_\theta|$ apertures are preallocated to each angle, then 
within each angle
only $|\mathcal{A}_\theta|!$ permutations are possible.
The total permutations $P_{p} = \prod_{\theta \in \Theta} |\mathcal{A}_\theta|! $, where $\sum_{\theta \in \Theta}|\mathcal{A}_\theta| = |\mathcal{A}|$, so unless all apertures are assigned to the same angle, $P_{p} < P_{d}$. %\Halmos
\endproof
Lemma \ref{lem:symGrowth} shows that there is redundancy in both models due to symmetry, but decision-based models grow faster than preallocated models for any plan with more than one active angle.
To mitigate the effect of symmetry in either model, it is desirable to include symmetry reduction techniques. Proposition \ref{preAppSym} shows that there are effective constraints for shrinking the decision space of both  preallocated and decision-based models. 
For the ease of notation, in what follows 
$P_{\text{d}}^{\text{(x)}}$ and  $P_{\text{p}}^{\text{(x)}}$ will denote the number of aperture permutations with the addition of symmetry-breaking constraint number (x), for the decision-based and preallocated models, respectively.

\begin{proposition} \label{preAppSym}
Using constraints \eqref{symCon1} and \eqref{symCon2} below, 
$P_d^{\eqref{symCon1}} = P_p^{\eqref{symCon1}} = P_p^{\eqref{symCon2}} = 1$ under the following conditions.

\noindent {\normalfont{a.}} 
$P_d^{\eqref{symCon1}} = P_p^{\eqref{symCon1}} = 1$ when $\sum_{\theta \in \Theta}\sum_{b \in \mathcal{B}_{\theta}} w_{b,a}$ is unique $\forall a \in \mathcal{A}$, and constraints \eqref{symCon1} are
\begin{align}
&\sum_{\theta \in \Theta} \sum_{b \in \mathcal{B}_{\theta}} w_{b,a}  \geq \sum_{\theta \in \Theta}\sum_{b \in \mathcal{B}_{\theta}} w_{b,a+1} ~~ \forall a \in \{1,\dots,|\mathcal{A}|-1\}. \label{symCon1}
\end{align}
{\normalfont{b.}} 
$P_p^{\eqref{symCon2}} = 1$ when  $\sum_{b \in \mathcal{B}_{\theta}} w_{b,a}$ is unique $\forall a \in \mathcal{A}$  within each angle  $\theta \in \Theta$. Constraints \eqref{symCon2} are 
\begin{align}
&\sum_{b \in \mathcal{B}_\theta} w_{b,a}  \geq \sum_{b \in \mathcal{B}_\theta} w_{b,a+1} ~~ \forall a \in \{1,\dots,|\mathcal{A}_\theta|-1\}, ~ \theta \in \Theta.  \label{symCon2}
\end{align}
There is no $P_d^{\eqref{symCon2}}$, as constraints \eqref{symCon2} do not apply to decision-based models.
\end{proposition}

\proof%{Proof.}
We discuss each condition in turn.

a. Constraints~\eqref{symCon1} sequence
total aperture fluence in a  monotonically decreasing order. Since only sequencing is considered, no unique, feasible plans are excluded from the decision space. If the uniqueness clause is violated, there will be
$P_d^{\eqref{symCon1}}=P_p^{\eqref{symCon1}}= \prod_{\eta \in \mathcal{F}} n_\eta!$ plan permutations, where $\mathcal{F}$ indexes the set of unique aperture fluence summations, $\sum_{b \in \mathcal{B}} w_{b,a}$, within a plan and $n_\eta$ is the total number of appearances of summation $\eta \in \mathcal{F}$.

b. 
In the case of preallocated models, tighter constraints~\eqref{symCon2} hold, as aperture intensities are only sorted within a set of angles. This further reduces $P_p^{\eqref{symCon2}}$ in the case of identical aperture sums, since if $\sum_{b \in \mathcal{B}_{\theta_1}} w_{b,a} = \sum_{b \in \mathcal{B}_{\theta_2}} w_{b,a}$, but
$\theta_1 \neq \theta_2$, the two sums are no longer included in the same index $\eta$. 

Constraints~\eqref{symCon2} can be shown to be incompatible with decision-based models
by contradiction: suppose constraints~\eqref{symCon2} do not impact the solution space of the decision-based model. Given a problem instance with $|\Theta| = 2$ beam angles, suppose there is a non-zero fluence value in aperture 1, angle 1, i.e., $\sum_{b \in \mathcal{B}_1} w_{b,1} > 0$. This, in turn, forces $\sum_{b \in \mathcal{B}_2} w_{b,1} = 0$. Every subsequent aperture must then be assigned a fluence of zero in angle 2 since the total fluence of angle 2 must not exceed the first aperture's total fluence value. The same issue occurs when the angles are reversed, ruling out any plan with radiation delivered from both angles, which greatly reduces the decision space, and contradicts our assumption. %\Halmos
\endproof

Proposition \ref{preAppSym} shows that tighter symmetry reduction constraints~\eqref{symCon2} can be applied only when preallocated models are used. Proposition \ref{propSymmetryBreak} shows that decision-based models can achieve the strong symmetry reduction seen in 
Proposition~\ref{preAppSym}b, using constraints~\eqref{Cons:SymBreak2} in the special case of two-angle treatments. Two-angle treatments are used in practice for certain regions of the body, such as tangential breast cancer RT \citep{purdie2011automated,purdie2014automation}. 

\begin{proposition}
\label{propSymmetryBreak}
When $|\Theta| = 2$, the following constraints result in a $P_d^{\eqref{Cons:SymBreak2}} = 1$ under the same conditions as $P_p^{\eqref{symCon2}} = 1$ in Proposition~\ref{preAppSym}b.
\begin{subequations}
\label{Cons:SymBreak2}
\begin{align}
\sum_{b \in \mathcal{B}_1} w_{b,a}  \geq \sum_{b \in \mathcal{B}_1} w_{b,a+1} && \forall a \in {1,\dots,|\mathcal{A}|-1}, \label{ang1Sym} \\
\sum_{b \in \mathcal{B}_2} w_{b,a}  \leq \sum_{b \in \mathcal{B}_2} w_{b,a+1} && \forall a \in {1,\dots,|\mathcal{A}|-1}. \label{ang2Sym}
\end{align}
\end{subequations}
\end{proposition}
\proof%{Proof.}
Each aperture $a$, has two associated total fluence values, $a_1 = \sum_{b \in \mathcal{B}_1} w_{b,a}$, and $a_2 = \sum_{b \in \mathcal{B}_2} w_{b,a}$. By definition, both $a_1$ and $a_2 \geq 0$. At least one of those equalities is strict by Constraints~\eqref{Cons:fullApt}.
By constraints~\eqref{ang2Sym}, any aperture with $a_2 = 0$ will be grouped and sequenced first, and any aperture with $a_1 = 0$ will be assigned to a group that is partitioned second by \eqref{ang1Sym}. Assuming no redundant summations, there is only one sequencing for the $a_1 > 0$ group, which also extends to the $a_2 > 0$ group. Any redundant summation values lead to possible plan permutations, but only within a given angle, as in the preallocated case in under constraints~\eqref{symCon2}.  
%\Halmos
\endproof

\begin{figure}[htb]
    \centering
    \includegraphics[width=0.7\textwidth]{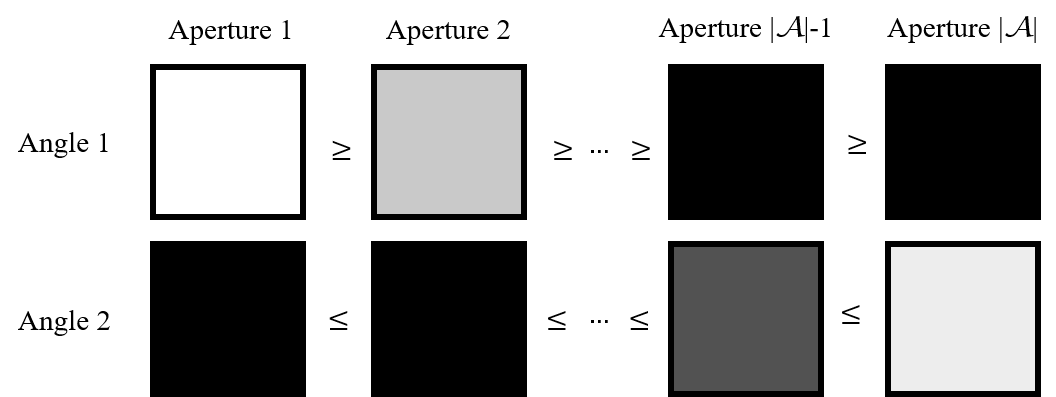}
    \caption{Visualization of intensity sorting by angle. Lighter boxes represent higher total beam intensities.}
    \label{fig:sortVis}
\end{figure}

Figure~\ref{fig:sortVis} provides the visual intuition behind Proposition \ref{propSymmetryBreak}. By grouping the darkened, or ``off'' angles, the apertures are partitioned into two groups. Within these groups, the sequencing of apertures does not impact the other group and therefore leads to an independent symmetry elimination at each angle. Incorporating constraints~\eqref{Cons:SymBreak2} into the decision-based RDAO model leads to a reduced solution space that eliminates identical solutions, potentially leading to faster solution times. For larger cases with more angles, however, only the less restrictive constraints~\eqref{symCon1} may be applied.

\subsubsection{Island Removal.}
\label{Sec:Island}
\begin{figure} [htb]
    \centering
    \begin{subfigure}[b]{0.5\textwidth}
        \centering
        \includegraphics[width=0.45\textwidth]{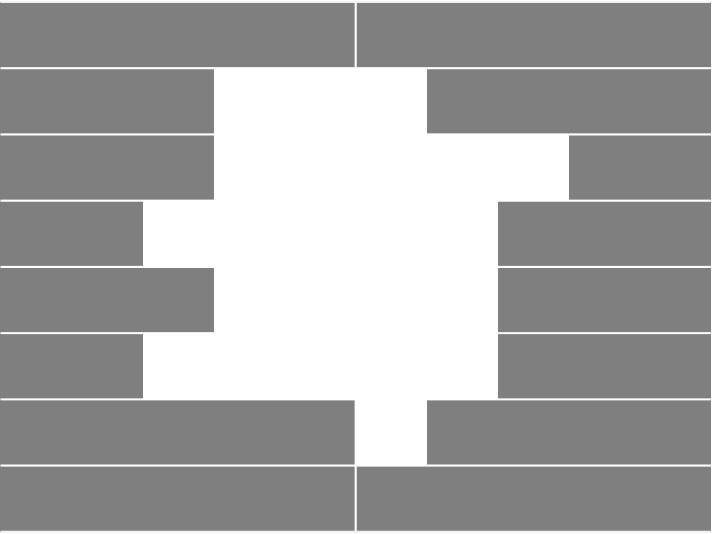}
        \caption{MLC leaves form a deliverable aperture}
        \label{legalIslands}
    \end{subfigure}%
    ~ 
    \begin{subfigure}[b]{0.5\textwidth}
        \centering
        \includegraphics[width=0.45\textwidth]{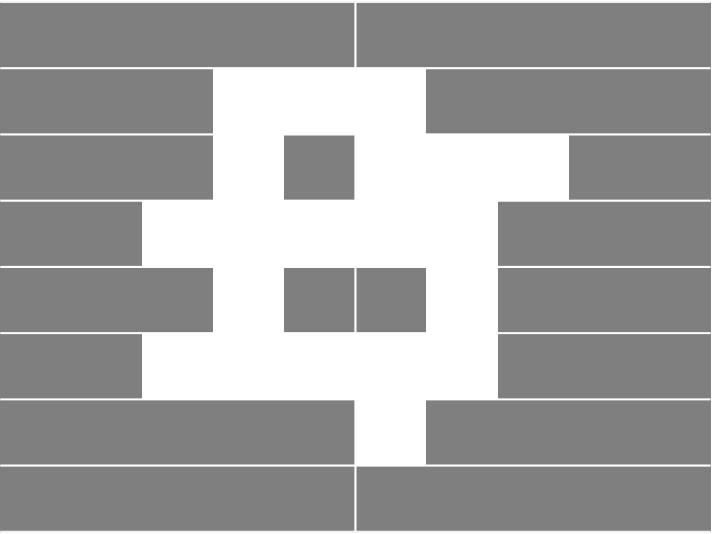}
        \caption{Undeliverable MLC setup with islands}
        \label{fig:illegalIsland}
    \end{subfigure}
    \caption{Possible MLC realizations with the current DAO constraints}
    \label{fig:islands}
\end{figure}

In practice, the deliverability of each aperture is limited by the linear beam modulation, performed by the MLC leaves. As such, the model must output beam setups that may be physically realized by the MLC. Mathematically, this means preventing breaks in the leaves, which result in detached sections or \textit{islands}, as shown in Figure~\ref{fig:islands}. 
In order to add these aperture requirements, the set of beamlets $\mathcal{B}$ is further partitioned from angles, into angle-dependent row and column coordinates, $\mathcal{Q}_\theta$ and $\mathcal{K}_\theta$, respectively. For the sake of exposition, it is assumed that the beam is an identical $|\mathcal{Q}|$ row $\times~|\mathcal{K}|$ column rectangle at each angle and that the leaves may traverse the entire span of the collimator from both the left and right sides. Moving forward, variables will be denoted using a dimension for each column, row and angle, instead of relying on the general beamlet identifier. The original beamlet identifier
% from the mapping 
$\mathcal{B} \leftarrow \lbrack \mathcal{Q}, \mathcal{K}, \Theta \rbrack$ can always be obtained using transformation
$b = \sum_{\theta' = 0}^{\theta-1}|\mathcal{B}_{\theta'}| + |\mathcal{K}|\times(q-1)+k$, where $|\mathcal{B}_{0} |= 0 $.

To restrict the leaves, the binary on-off constraints $x_{q,k,\theta,a}$, used to enforce uniformity, can be joined by two additional sets of binary variables, $l_{q,k,\theta,a}$ and $r_{q,k,\theta,a}$. These variables represent the continuous leaves extended from the left and right side of the collimator, respectively, and impact $x_{q,k,\theta,a}$ as follows: 
\[  x_{q,k,\theta,a}  = 
 \left\{
\begin{array}{ll}
      1 & \text{if } l_{q,k,\theta,a} =1 \text{ and } r_{q,k,\theta,a} = 1,  \\
      0 & \text{otherwise}. \\
\end{array} 
\right. \]
A left leaf is open (i.e., not extended over a beamlet) when $l = 1$, and similarly, an open right leaf is indicated by $r=1$. If $l = 1$ and $r = 1$ a beamlet is open, meaning it is on, or $x =1$. If either $l$ or $r = 0$, a beamlet is closed. Both left and right leaves cannot cover the same beamlet simultaneously, meaning they cannot both be 0, or $l + r \geq 1$. 
The following constraints create continuous non-overlapping leaves.
\begin{subequations}
\label{Cons:islands}
\begin{align}
&l_{q,k+1,\theta,a} \geq l_{q,k,\theta,a} && \forall k \in \mathcal{K}', q \in \mathcal{Q}, \theta \in \Theta, a \in \mathcal{A}, \label{leftFinger} \\
&r_{q,k,\theta,a} \geq r_{q,k+1,\theta,a} && \forall k \in \mathcal{K}', q \in \mathcal{Q}, \theta \in \Theta, a \in \mathcal{A}, \label{rightFinger} \\
&x_{q,k,\theta,a} = -1+ l_{q,k,\theta,a}+r_{q,k,\theta,a} && \forall k \in \mathcal{K}, q \in \mathcal{Q}, \theta \in \Theta, a \in \mathcal{A}, \label{eitherOr}\\ 
&l_{q,k,\theta,a},r_{q,k,\theta,a} \in \{0,1\} && \forall k \in \mathcal{K}, q \in \mathcal{Q}, \theta \in \Theta, a \in \mathcal{A}, \label{islandBin}
\end{align}
\end{subequations}
where, $\mathcal{K}' = \{1, ..., |\mathcal{K}| -1\}$.

An alternative approach to island removal, which uses a larger number of variables but fewer constraints, was proposed for FMO leaf sequencing by \citep{boland2004minimizing}. While their constraints result in the same feasible set of $\textbf{x}$ and $\textbf{w}$ variables, we found that our proposed method runs faster in the full model empirically, possibly due to the reduced number of binary variables. 

\subsubsection{Vertical Continuity.}
\label{sec:VertCont}
So far, DAO constraints have been derived to ensure the deliverability of the plans for an MLC. 
There are, however, certain features that clinicians might prefer to see in a final plan, that are not strict mechanical deliverability requirements. For instance, multiple separated groups of beamlets within an aperture can be undesirable due to a possible leakage between adjacent closed leaves. This phenomenon can be largely avoided through the creation of continuous segments that are mostly covered by jaws. 

\begin{figure} [h!]
    \centering
    \begin{subfigure}[b]{0.33\textwidth}
        \centering
        \includegraphics[width=0.65\textwidth]{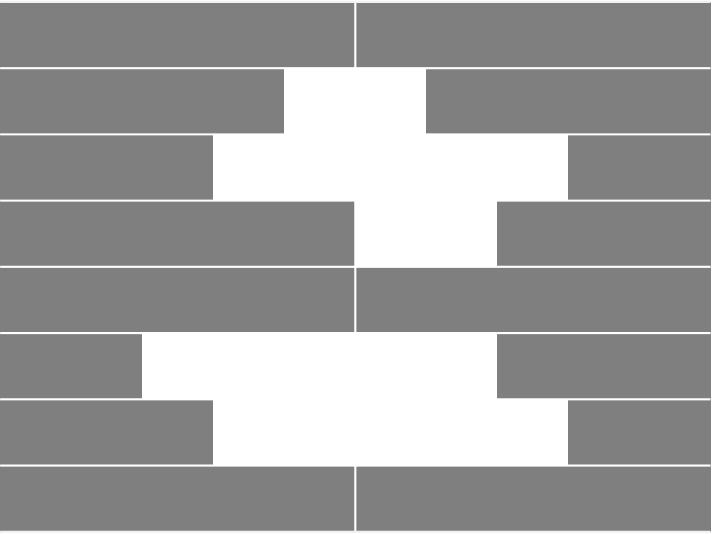}
        \caption{Vertically separated}
        \label{VertSep}
    \end{subfigure}%
    ~  
    \begin{subfigure}[b]{0.33\textwidth}
        \centering
        \includegraphics[width=0.65\textwidth]{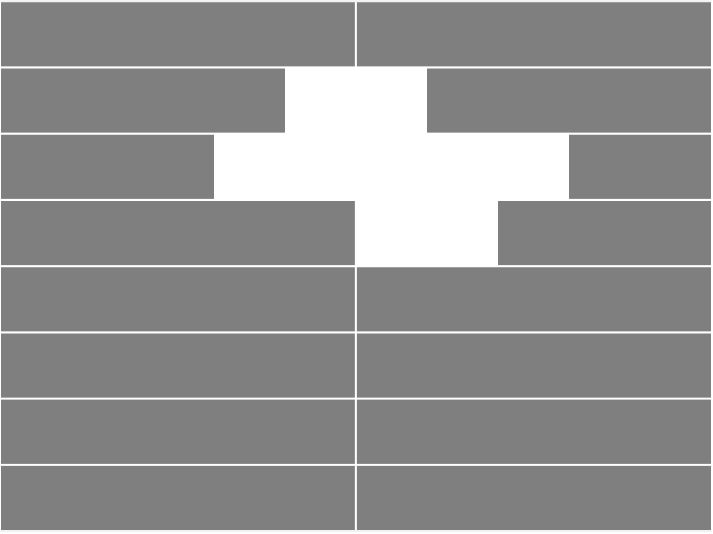}
        \caption{First aperture}
        \label{VertSep_split1}
    \end{subfigure}%
    ~ 
    \begin{subfigure}[b]{0.33\textwidth}
        \centering
        \includegraphics[width=0.65\textwidth]{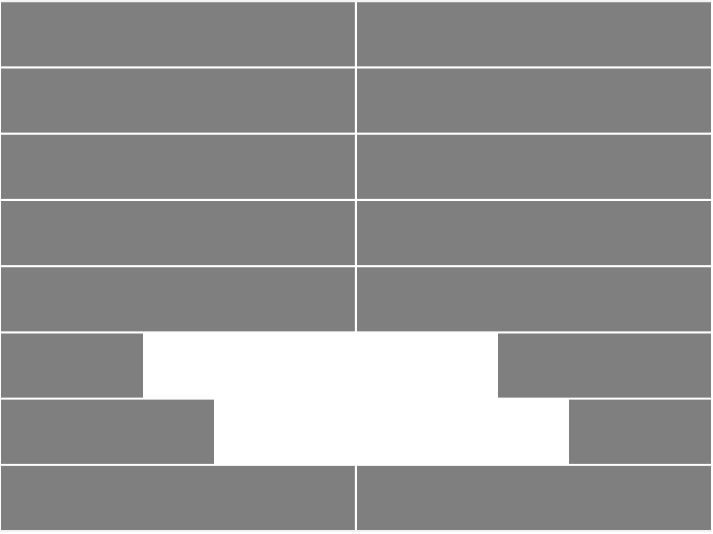}
        \caption{Second aperture}
        \label{VertSep_split2}
    \end{subfigure}  
      
	\par\bigskip

    \begin{subfigure}[b]{0.33\textwidth}
        \centering
        \includegraphics[width=0.65\textwidth]{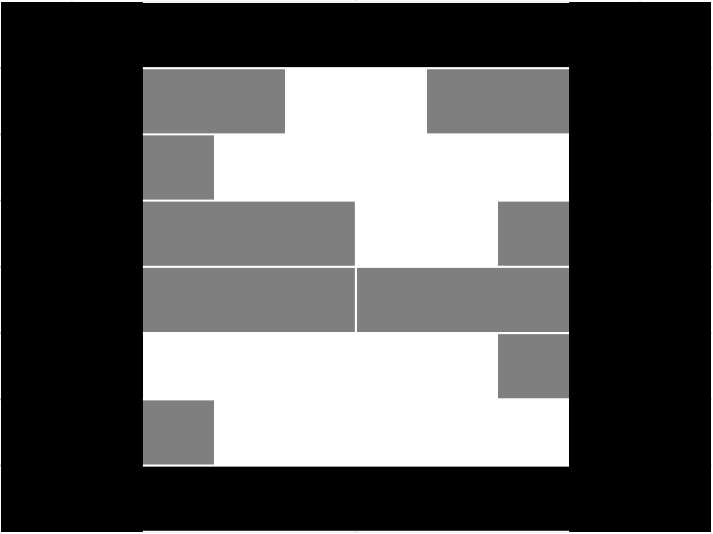} 
        \caption{Original jaw coverage}
       \label{VertSep_jaws}
    \end{subfigure}%
    ~ 
    \begin{subfigure}[b]{0.33\textwidth}
        \centering
        \includegraphics[width=0.65\textwidth]{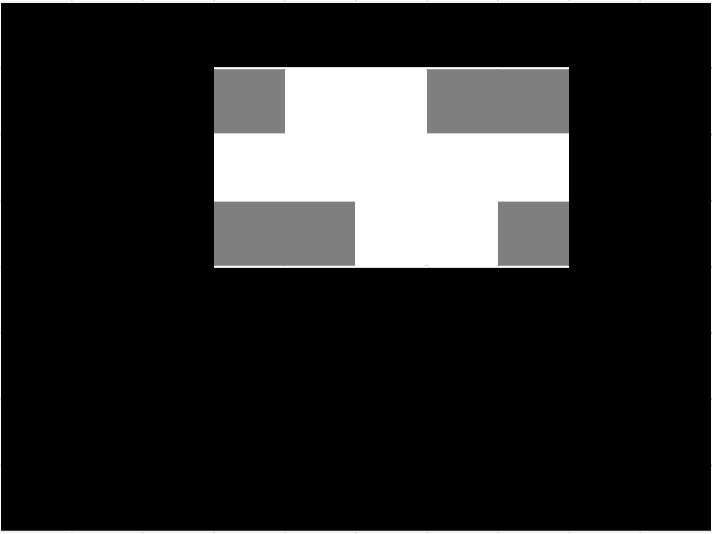}
        \caption{First jaw coverage}
           \label{legalsplit_jaws1}
    \end{subfigure}%
    ~ 
    \begin{subfigure}[b]{0.33\textwidth}
        \centering
        \includegraphics[width=0.65\textwidth]{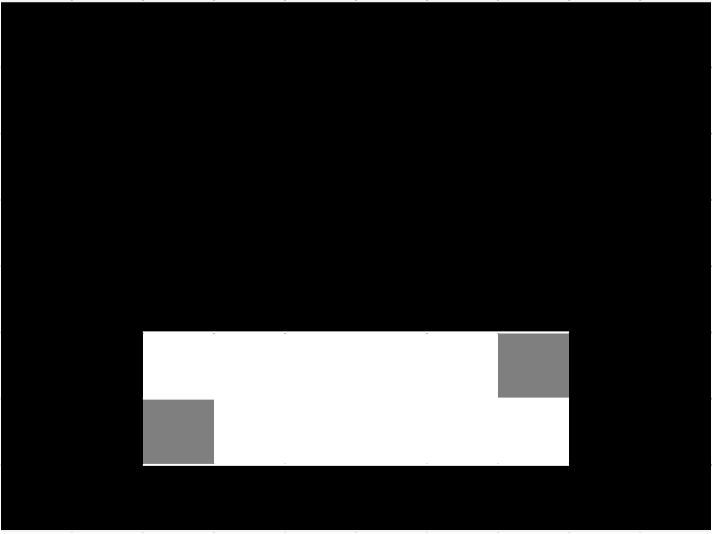}
        \caption{Second jaw coverage}
           \label{legalsplit_jaws2}        
    \end{subfigure}    
    \caption{Deliverable aperture with clinically undesirable vertical break, separated to two apertures. Images (d-f) use darker color to show jaw coverage.}
    \label{multipleSegs}
\end{figure}

A clinician presented with the split aperture shown in Figure~\ref{VertSep}, for example, would typically, manually create two separate apertures, as in Figures~\ref{VertSep_split1} and \ref{VertSep_split2}. 
To see why, contrast the leakage protection by the jaws in Figure~\ref{VertSep_jaws}, with that of the combined Figures~\ref{legalsplit_jaws1} and~\ref{legalsplit_jaws2}.
This change, however, would increase the total number of apertures; a parameter previously controlled for in the model. 

Constraints~\eqref{Cons:vert} enforce the same principle by disallowing vertically-disconnected apertures, through the introduction of $2\times |\mathcal{Q}| \times |\Theta| \times |\mathcal{A}|$ new binary variables. Variables $\bar{\textbf{j}}$ and $\underline{\textbf{j}}$ restrict upper and lower jaw motion, respectively. 
Variable $\bar{j}_{q,\theta,a}$ indicates if the upper jaw at angle $\theta$, aperture $a$ is blocking row $q$ (i.e., $\bar{j}_{q,\theta,a}=0$) or open (i.e., $\bar{j}_{q,\theta,a}=1$). Similarly,  $\underline{j}_{q,\theta,a}$ does the same for the lower jaw. Summary variable $j_{q,\theta,a}$ indicates whether or not a row is active. Together, these constraints prevent vertical breaks between active rows. 
\begin{subequations}
\label{Cons:vert}
\begin{align}
& j_{q,\theta,a} =  -1 + \bar{j}_{q,\theta,a} + \underline{j}_{q,\theta,a} && \forall q \in \mathcal{Q}, \theta \in \Theta, a \in \mathcal{A}, \label{sumVar} \\
& j_{q,\theta,a} \leq \sum_{k \in \mathcal{K}} x_{q,k,\theta,a} &&  \forall q \in \mathcal{Q}, \theta \in \Theta, a \in \mathcal{A}, \label{singleBeamlet} \\
& |\mathcal{K}| \times  j_{q,\theta,a} \geq \sum_{k \in \mathcal{K}} x_{q,k,\theta,a} &&  \forall q \in \mathcal{Q}, \theta \in \Theta, a \in \mathcal{A}, \label{maxWholeRow} \\
& \bar{j}_{q,\theta,a} \leq \bar{j}_{q+1,\theta,a}  && \forall q \in \mathcal{Q}', \theta \in \Theta, a \in \mathcal{A}, \label{upperVert}\\ 
& \underline{j}_{q+1,\theta,a} \leq \underline{j}_{q,\theta,a}  && \forall q \in \mathcal{Q}', \theta \in \Theta, a \in \mathcal{A}, \label{lowerVert}\\ 
& j_{q,\theta,a}, \bar{j}_{q,\theta,a},\underline{j}_{q,\theta,a} \in \{0,1\} && \forall q \in \mathcal{Q}, \theta \in \Theta, a \in \mathcal{A}, \label{vertBin}
\end{align}
\end{subequations}
where $\mathcal{Q}' = \{1, ..., |\mathcal{Q}| -1\}$.

\subsubsection{Horizontal Continuity.}
\label{sec:HorCont}
\begin{figure} [htb]
    \centering
    \begin{subfigure}[b]{0.5\textwidth}
        \centering
        \colorbox{white}{\includegraphics[width=0.6\textwidth]{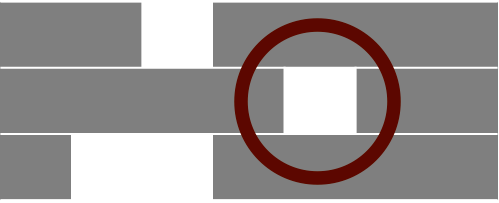}}
        \caption{Leaf collision}
        \label{collision}
    \end{subfigure}%
    ~ 
    \begin{subfigure}[b]{0.5\textwidth}
        \centering
        \includegraphics[width=0.6\textwidth]{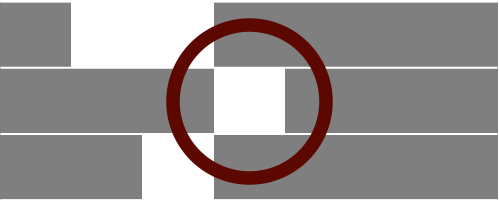}
        \caption{Disconnected rows}
        \label{singleton}
    \end{subfigure}
    \caption{Examples of undesirable behavior between pairs of adjacent rows}
    \label{leafPositions}
\end{figure}

In the same vein as vertical continuity, there is a clinical incentive to avoid difficult, leakage-prone setups such as leaf collisions and disconnected rows. Leaf collisions can occur when the right leaves extend beyond adjacent row left leaves (or vice-versa) as shown in Figure~\ref{collision}. If the leaves become even slightly misaligned, they might collide during setup, so given a choice, planners typically prefer to mitigate this risk.

Disconnected rows have no vertical beamlet connection between adjacent rows, leading to horizontally-separated segments, as shown in Figure~\ref{singleton}. Since leakage primarily arises around the edges of leaves, plans with very small and separated sections can lead to unwanted leakage. 

While previous work on leaf sequencing has introduced preventative measures for collisions, their models still allow rows to become disconnected \citep{boland2004minimizing}. To prevent all horizontal detachment, the following set of constraints is proposed,
\begin{subequations}
\begin{align}
j_{q,\theta,a} + j_{q-1,\theta,a} - \sum_{\delta=k+1}^{|\mathcal{K}|} x_{q,\delta,\theta,a} \leq 1 + \sum_{\delta=1}^{k} x_{q-1,\delta,\theta,a} && \forall k \in \mathcal{K}, q \in \mathcal{Q}'', \theta \in \Theta, a \in \mathcal{A}, \label{noSingletonsLeft} \\
j_{q,\theta,a} + j_{q-1,\theta,a} - \sum_{\delta=1}^{|\mathcal{K}|-k} x_{q,\delta,\theta,a} \leq 1 + \sum_{\delta=|\mathcal{K}|-k+1}^{|\mathcal{K}|} x_{q-1,\delta,\theta, a} && \forall k \in \mathcal{K}, q \in \mathcal{Q}'', \theta \in \Theta, a \in \mathcal{A}, \label{noSingletonsRight} 
\end{align}
\end{subequations}
where $\mathcal{Q}'' = \{2,\dots,|\mathcal{Q}|\}$. Constraints~\eqref{noSingletonsLeft} and \eqref{noSingletonsRight} ensure from the left and right side, respectively, that all active, adjacent rows share at least one active beamlet, before reaching the total number of active beamlets in that row. 
The additional $j_{q,\theta,a}$ terms represent edge cases, relaxing the constraints when one or both of the rows in a pair are off. All other edge cases are eliminated by the no islands constraints~\eqref{Cons:islands}.

\subsection{Complete RDAO Model}
\label{sec:abridgedMod}
Putting all the constraints in this section together, the full, proposed robust direct aperture optimization (RDAO) model can be derived as follows,
\begin{align}
\tag{\textbf{RDAO}} \label{RDAOMod}\\
\text{min~} & \sum_{s\in \{\cT,\cH\}}\frac{c_s}{|\mathcal{V}_s|}\sum_{v \in \mathcal{V}_s} \sum_{b \in \mathcal{B}} \sum_{i \in \mathcal{I}} \sum_{a \in \mathcal{A}} p_i D_{v,b,i} w_{b,a}  \nonumber \\ 
\text{s.t.~}& 
\eqref{RobModInf}, \eqref{RobPos} \text{ with substitution } \eqref{sumW}, \tag{Robustness} \\
& \eqref{blockedBeam},\eqref{openBeam},\eqref{openBeam2},\eqref{PositiveDose}, \eqref{binCon}, \tag{Uniformity} \\
& \eqref{angleChoice},\eqref{oneAngle}  \text{ and optionally } 
\eqref{symCon1}  \text{ or }  \eqref{ang1Sym}, \eqref{ang2Sym},
\tag{Aperture Selection} \\
& \eqref{leftFinger}, \eqref{rightFinger}, \eqref{eitherOr},\eqref{islandBin}, \tag{Island Removal}  \\
& \eqref{sumVar},\eqref{singleBeamlet},\eqref{maxWholeRow},\eqref{upperVert}, \eqref{lowerVert},\eqref{vertBin}, \tag{Vertical Continuity} \\
& \eqref{noSingletonsLeft},\eqref{noSingletonsRight}. \tag{Horizontal Continuity} 
\end{align}
This same \eqref{RDAOMod} model with explicit constraints is detailed in 
%Online 
Appendix 
%B.
\ref{Sec:fullMod}.

%%%%%%%%%%%%%%%~~~~~  SECTION 3 ~~~~~%%%%%%%%%%%%%%%

\section{Candidate Plan Generation Heuristic}
\label{Sec:warmStart}

The RDAO model is an extremely large-scale MIP model, which makes it computationally challenging to solve. Even finding an initial feasible solution can be non-trivial, often resulting in long search times and very large initial optimality gaps. For this reason, we propose the candidate plan generation (CPG) heuristic as a method for finding high-quality, deliverable plans. The CPG heuristic is similar to some of the smoothing techniques in the literature such as the nearest neighbor-pair smoothing in \cite{saberian2017spatiotemporally}, which was inspired by the objective smoothing term in \cite{breedveld2006fast} 
and earlier work by
\cite{webb1998inverse}. 
Unlike the CPG method, however, the \cite{saberian2017spatiotemporally} nearest neighbor constraints prevent the sharp dosage drops necessary for robust optimization, while the objective smoothing term in \cite{breedveld2006fast} results in a loss of control over the final
number of apertures.

The CPG plans generated by the heuristic can either be used 
as a standalone approach,
or translated to \eqref{RDAOMod} variables, providing an initial incumbent solution to warm start the solver. Section~\ref{Sec:wsS1} details the three steps of the heuristic, while Section~\ref{Sec:warmStartBridge} discusses the process of translating CPG output into a feasible incumbent solution for \eqref{RDAOMod}.

\subsection{The Three-Step CPG Heuristic} 
\label{Sec:wsS1}
\begin{figure}[htbp]
    \centering
    \includegraphics[width=0.6\textwidth]{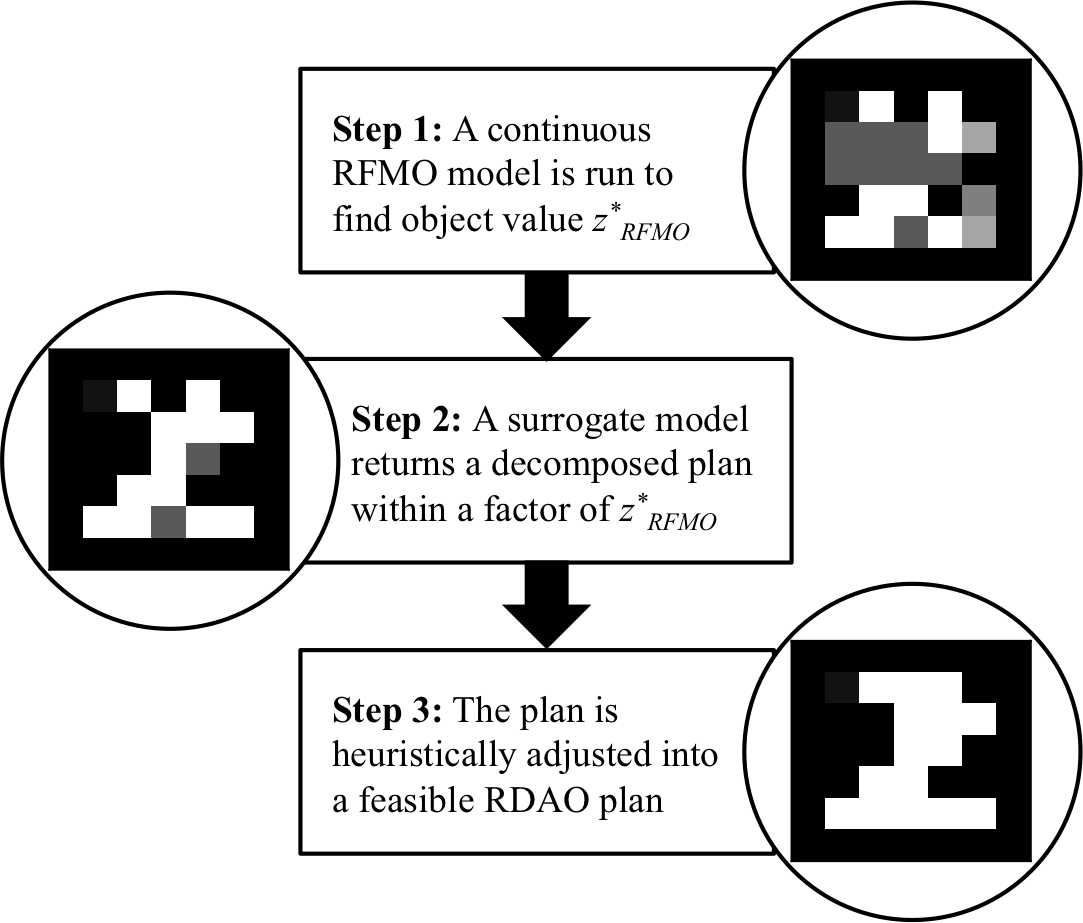}
    \caption{Outline of the 3-step candidate plan generation heuristic, with sample output aperture depictions.}
    \label{fig:warmStartFlow}
\end{figure}

The CPG heuristic, outlined in Figure~\ref{fig:warmStartFlow}, is a three-step process in which each step is a continuous LP or a rapid algorithm. 
The process for each step is as follows.

\noindent 
\textbf{Step 1: Find a Lower Bound}
by running the continuous \eqref{RobMod}.

\noindent \textbf{Output:} 
$z^*_{RFMO}$ -- An optimal \eqref{RobMod} objective function value.

\noindent 
\textbf{Step 2: Run a Surrogate Model} 
\label{Sec:wsS2}
to find a lower bound on each of the CPG fluence intensities, $\underline{\textbf{w}}^{cpg}$, within some factor of $z^*_{RFMO}$.
This surrogate model approximates the binary segment-uniformity constraints with a min-max objective, which has a leveling effect on every active beamlet in a given aperture. 
\begin{subequations}
\begin{align}
\tag{\textbf{CPG-S}} \label{WSMod}\\
\min~ & \alpha \sum_{\theta \in \Theta} \sum_{a \in \mathcal{A'}} m_{\theta,a} + (1-\alpha) \sum_{s\in \{\cT,\cH\}}\frac{c_s}{|\mathcal{V}_s|}\sum_{v \in \mathcal{V}_s} \sum_{b \in \mathcal{B}} \sum_{i \in \mathcal{I}} \sum_{a \in \mathcal{A'}} p_i D_{v,b,i} \underline{w}_{b,a}^{cpg} \label{newWarmObj} \\
\text{s.t.~} & \sum_{b \in \mathcal{B}} \sum_{i \in \mathcal{I}}  \sum_{a \in \mathcal{A'}}  \tilde{p}_i D_{b,v,i} \underline{w}_{b,a}^{cpg}  \geq L_{v}  \hspace{132pt} \forall v \in \mathcal{V}_T, ~ \forall \tilde{\textbf{p}} \in \mathcal{P}, \label{WarmTum} \\ 
& m_{\theta,a} \geq \underline{w}_{b,a}^{cpg} \hspace{220pt}  \forall  b \in \mathcal{B_\theta}, a \in \mathcal{A'}, \theta \in \Theta, \label{warmMinMax} \\
& \underline{w}_{b,a}^{cpg} \geq 0  \hspace{236pt}  \forall b \in \mathcal{B}, a \in \mathcal{A},  \label{warmPos}
\end{align}
\end{subequations}
where $\mathcal{A'} =   \{1,\dots,\frac{|\mathcal{A}|}{|\Theta|}\}$, $\mathcal{A}$ is the desired set of apertures in the final model, 
$0\leq\alpha\leq1$ is a tuneable parameter, and $m_{\theta, a}$ is the maximum intensity value at angle $\theta$ in aperture $a$. 
The model resembles the continuous \eqref{RobMod}, with the following three differences: First, the new objective function \eqref{newWarmObj} has two terms; the first minimizes the sum of maximum dose in each aperture (in use when $\alpha > 0$), while the second is identical to the \eqref{RobMod} model's objective function \eqref{RobObj} (in use when $\alpha < 1$). The weighting term $\alpha$ dictates the balance between uniformity and original objective function minimization. 
Second, the min-max property is enforced in constraints~\eqref{warmMinMax}. Finally, the DAO variables $\underline{w}_{b,a}^{cpg}$ must be used, rather than $\omega_b$, to store intensities at each of the $|\mathcal{A'}|$ apertures, which are implicitly preallocated equally to angles. This assumption can be loosened with additional decision variables, but doing so within the CPG is outside the scope of this paper.

\noindent \textbf{Output:} $\underline{w}_{b,a}^{cpg}$ -- A lower bound on a set of deliverable intensities, and $m_{\theta, a}$ -- The maximum intensity at each active aperture-angle pair.

\noindent 
\textbf{Step 3: Gap Filling},
\label{Sec:wsS3}
as outlined in Algorithm 1, is done to enforce deliverability constraints while remaining as close as possible to the Step~2 optimal plan. The goal is to output a set of intensities appropriate for both clinical use and warm-starting a solver.  

\SetKwBlock{Repeat}{Optional}{}

\newcounter{algoline}
\newcommand\Numberline{\refstepcounter{algoline}\nlset{\thealgoline}}
\AtBeginEnvironment{algorithm}{\setcounter{algoline}{0}}
\IncMargin{2em}
\begin{algorithm}[h]
\caption{Gap Filling}\label{gapFill}
\nlset{A1.1} \label{A1s1}%
\textbf{initialize} intensities $w_{b,a}^{cpg}$ to 0 for all $b \in \mathcal{B}$ and $a\in \mathcal{A}$.\\
\textbf{create} $s_a$: a symmetry-breaking sequencing of the $|\mathcal{A}'|\times|\Theta|$ apertures in $m_{\theta, a}$.\\
 \For{aperture $a \in \mathcal{A}'$ and angle $\theta \in \Theta$}{ 
  \nlset{A1.2} \label{A1s2} %\Numberline 
  \For{each row $q$}{ 
  \textbf{initialize} first and last active column $\underline{k}_{q,\theta,a} \leftarrow \{\}$ and $\bar{k}_{q,\theta,a} \leftarrow \{\}$ \\
   \textbf{if }$\exists~ \underline{w}_{q,k,\theta,a}^{cpg}>0$ \textbf{then} $\underline{k}_{q,\theta,a} \leftarrow \min_k\{\underline{w}_{q,k,\theta,a}^{cpg}>0\}$, $\bar{k}_{q,\theta,a} \leftarrow \max_k\{\underline{w}_{q,k,\theta,a}^{cpg}>0\}.$
   }
\nlset{A1.3} \label{A1s3} 
\textit{Optional:} \textbf{initialize} $q' \leftarrow 0$, $\underline{k}' \leftarrow 1$, $\bar{k}' \leftarrow |\mathcal{K}|$ to track the last active row.\\ 
\For{each active ($\underline{k}_{q,\theta,a} \neq \varnothing$) row $q$}{
\textbf{if} row $q$ starts after the last ($q'$) finishes, $\underline{k}_{q,\theta,a} > \bar{k}'  $ \textbf{then} $\underline{k}_{q,\theta,a} \leftarrow \bar{k}'.$ \\
\textbf{if} row $q$ finishes before $q'$ starts, $\bar{k}_{q,\theta,a} < \underline{k}' $ \textbf{then} $\bar{k}_{q,\theta,a} \leftarrow \underline{k}'.$\\
 \If{($q>0$ \textbf{and} $q-q'>1$)}{ %\textbf{do} 
 $\bar{k}_{t,\theta,a} \leftarrow \underline{k}_{t,\theta,a} \leftarrow \max{\{\underline{k}',\underline{k}_{q,\theta,a}\}} ~~~ \forall t \in \{q'+1, \dots, q-1\}.$}
\textbf{update} $q' \leftarrow q$,  $\underline{k}' \leftarrow \underline{k}_{q,\theta,a}$ and $\bar{k}' \leftarrow \bar{k}_{q,\theta,a}$.
}
\nlset{A1.4} \label{A1s4} 
\textbf{for} each row $q$ 
\textbf{assign} $w^{cpg}_{q,k,\theta,s_a} \leftarrow m_{\theta,a} ~~~ \forall k \in \{\underline{k}_{q,\theta,a}, \dots, \bar{k}_{q,\theta,a}\}.$
} 
\end{algorithm}
\DecMargin{2em}
The gap filling algorithm has four main parts. The first part, \ref{A1s1}, is initialization, where fluence values are set to zero, and apertures are sequenced based on the desired symmetry-breaking pattern.  
\ref{A1s2} contains the search for the index of the first and last active column in each row (all middle indices must be active to avoid islands).
\ref{A1s3} is an optional segment, 
wherein row endpoints are adjusted to ensure that there are no gaps between adjacent rows, either horizontally or vertically. Finally, \ref{A1s4}, fills in uniform fluence values to the assigned aperture, based on each row's endpoints.

\noindent \textbf{Output:} $w^{cpg}_{b,a}$ - A deliverable direct aperture plan.

\noindent Proposition~\ref{warmProp} formally demonstrates that the solution generated by the CPG heuristic meets all the requirements of the~\eqref{RDAOMod} model, and is therefore feasible.

\begin{proposition}
\label{warmProp}
The solution of the CPG heuristic is feasible for \eqref{RDAOMod}.
\end{proposition}
\proof%{Proof.}  
It is sufficient to show that the new $w^{cpg}_{b,a}$  values generated in Algorithm 1 meet each of the following \eqref{RDAOMod} constraints:

\noindent \textbf{\textit{Robustness}}: Plan $\underline{w}_{b,a}^{cpg}$ from Step~2 satisfies robust dose constraints~\eqref{RobModInf} and \eqref{RobPos}, as they are also included in the surrogate model as constraints~\eqref{WarmTum} and \eqref{warmPos}. Since intensity values are increased in the heuristic, never decreased, we can guarantee that ${\sum_{a\in\mathcal{A}} w_{b,a}^{cpg} \geq \sum_{a\in\mathcal{A'}} \underline{w}_{b,a}^{cpg},} ~ \forall b \in \mathcal{B}$, and therefore $w_{b,a}^{cpg}$ must also satisfy robust constraints. 
\noindent \textbf{\textit{Uniformity}}: All nonzero values assigned to each aperture-angle pair in $w_{b,a}^{cpg}$ take on a value of $m_{\theta,a}$ in 
A1.4, hence uniformity constraints~\eqref{blockedBeam}--\eqref{binCon} are met.

\noindent \textbf{\textit{Aperture Selection}}: Since $|\mathcal{A'}| =  \frac{|\mathcal{A}|}{|\Theta|}$, the $|\mathcal{A'}|\times|\Theta|$ apertures from $\underline{w}_{b,a}^{cpg}$ are each assigned to one of the $|\mathcal{A'}|$ apertures in $w_{b,a}^{cpg}$ in A1.4, therefore $
\sum_{\theta \in \Theta} u_{a,\theta} = 1 ~~ \forall a \in \mathcal{A}$ is satisfied, and by extension, constraints~\eqref{oneAngle} and \eqref{angleChoice}.
Dose intensities are assigned in 
A1.4 according to sequencing $s_a$ introduced in 
A1.1, ensuring that any desired symmetry-breaking sequencing is also conserved, satisfying \eqref{symCon1} or \eqref{ang1Sym} and \eqref{ang2Sym}.

\noindent \textbf{\textit{Island Removal}}: Filling in all intensities between the first and last active columns  in each active row, $\underline{k}_{q,\theta,a}$ and $\bar{k}_{q,\theta,a}$, identified in 
A1.2 and enforced in 
A1.4, guarantees no islands can form, thereby satisfying constraints~\eqref{leftFinger}--\eqref{islandBin}.

\noindent \textbf{\textit{Vertical and Horizontal Continuity}}: When 
A1.3 is applied, horizontal continuity is guaranteed by adjusting non-overlapping row-endpoints $\underline{k}_{q,\theta,a}$ or $\bar{k}_{q,\theta,a}$ to overlap with the last active row by one beamlet, satisfying constraints \eqref{noSingletonsLeft} and \eqref{noSingletonsRight}. Rows separated vertically 
are then connected using a single-beamlet column, satisfying \eqref{sumVar}--\eqref{vertBin}.
We claim $\max{\{\underline{k}',\underline{k}_{q,\theta,a}\}}$ is a connecting column where horizontal continuity is guaranteed, which can be shown by contradiction. 
If the row with a further starting point does not share a common column with the earlier starting point, either $\underline{k}_{q,k,\theta,a} > \bar{k}'$ or $\underline{k}' > \bar{k}_{q,k,\theta,a}$, implying horizontal adjustments have not been made, contradicting the initial claim. %\Halmos
\endproof

\subsection{Full Warm Start Solution Generation}
\label{Sec:warmStartBridge}
In order to use deliverable intensities $w^{cpg}_{b,a}$ as a warm start for a commercial solver, the remaining mixed-integer \eqref{RDAOMod} variable values must be generated.
Since Proposition \ref{warmProp} has already shown that the plan intensities are feasible, the remaining variables can be 
made feasible by construction. 
Warm start generation algorithm, Algorithm 2, constructs a feasible starting point for the full \eqref{RDAOMod} model from deliverable intensities.
\IncMargin{2em}
\begin{algorithm}[h]
\caption{Warm Start Generation}\label{warmStA}
\nlset{A2.1} \label{A2s1}%
\textbf{initialize} $x_{b,a}^{cpg}$, $f_a^{cpg}$, $u_{a,\theta}^{cpg}$, $l_{q,k,\theta,a}^{cpg}$, $j_{q,\theta,a}^{cpg}$ and $\bar{j}_{q,\theta,a}^{cpg}$ to \textbf{0}.\\
\textbf{initialize}  $r_{q,k,\theta,a}^{cpg}$ and $\underline{j}_{q,\theta,a}^{cpg}$ to \textbf{1}. \\
\textbf{initialize} \textit{beamOn}, \textit{beamOff}, \textit{rowOn}, \textit{rowOff} to 0.\\
\For{all apertures $a \in \mathcal{A}$}{
\nlset{A2.2} \label{A2s2}
\textbf{set}  $f_a^{cpg} \leftarrow \max_a\{w^{cpg}_{b,a}\}$ and $\theta \leftarrow angle(\max_a\{w^{cpg}_{b,a}\})$,\\
\textbf{set}  $u_{a,\theta}^{cpg} \leftarrow 1$,\\
\textit{beamOn} $\leftarrow 0$, \textit{beamOff} $\leftarrow 0$.\\
\For{each row $q \in \mathcal{Q}$ in each angle $\theta \in \Theta$}{
 \textbf{set} \textit{rowOn} $\leftarrow 0$, \textit{rowOff} $\leftarrow 0$. \\
\nlset{A2.3} \label{A2s3}
\For{each columns $k \in \mathcal{K}$}{
\textbf{if} $w^{cpg}_{q,k,\theta,a} > 0$ \textbf{then} $x_{q,k,\theta,a}^{cpg} \leftarrow 1$ and \textit{rowOn} $\leftarrow 1$. \\
\textbf{else if} \textit{rowOn} \textbf{then} \textit{rowOff} $\leftarrow 1$.\\
\textbf{if} \textit{rowOn} \textbf{then} $l_{q,k,\theta,a}^{cpg} \leftarrow 1$.\\
\textbf{if} (\textit{rowOn} \textbf{\textit{and}} \textit{rowOff}) \textbf{then} $r_{q,k,\theta,a}^{cpg} \leftarrow 0$.
}
\nlset{A2.4} \label{A2s4}
\textit{Optional:} \textbf{if} \textit{rowOn} \textbf{then} $j_{q,\theta,a}^{cpg} \leftarrow 1$ and  \textit{beamOn} $\leftarrow 1$.\\
\textbf{else if} \textit{beamOn} \textbf{then} \textit{beamOff} $\leftarrow 1$.\\
\textbf{if} \textit{beamOn} \textbf{then} $\bar{j}_{q,\theta,a}^{cpg} \leftarrow 1$.\\
\textbf{if} (\textit{beamOn} \textbf{\textit{and}} \textit{beamOff}) \textbf{then} $\underline{j}_{q,\theta,a}^{cpg} \leftarrow 0$.\\
}}
\nlset{A2.5} \label{A2s5}
\textbf{run} \eqref{RobMod}  fixing $\omega_b = \sum_{a \in \mathcal{A}} w_{b,a}^{cpg}, ~ \forall b \in \mathcal{B}$, for robust counterpart variables.
\end{algorithm}
\DecMargin{2em}

The algorithm has five parts. In the first part, A2.1, all variables are initialized to 0, except the right leaves and lower jaw, which are initialized to 1, or the ``on'' or ``open'' position. 
Conceptually, this means that all left leaves and upper jaws are assumed to completely span the beam when the algorithm begins. In 
A2.2, active aperture angles and their corresponding fluence values are identified. Flag pairs \textit{beamOn} and \textit{beamOff}, and \textit{rowOn} and \textit{rowOff}, which are used detect the start and end of apertures and beamlet-openings, respectively, are also reset. 
A2.3 handles the beamlet activation binaries, including the deactivation of left leaves and the activation of right leaves in any row with active beamlets. 
A2.4 does the same with the upper and lower jaws, respectively, when vertical and horizontal continuity constraints have been activated.
Finally, since robust optimization is typically done using a robust counterpart, A2.5  
suggests running \eqref{RobMod}
with fixed values based on $w_{b,a}^{cpg}$, to find any dual variable values, yielding a complete candidate solution. 

\noindent \textbf{Output:} variable values for \eqref{RDAOMod}, including: $\{\textbf{x}, \textbf{f}, \textbf{u}, \textbf{l}, \textbf{r}, \textbf{j}, \bar{\textbf{j}}, \underline{\textbf{j}}\}$.

%%%%%%%%%%%%%%%~~~~~  SECTION 4 ~~~~~%%%%%%%%%%%%%%%

\section{Results and Discussion}
\label{Sec:Results}
The \eqref{RDAOMod} model and CPG heuristic were tested on five clinical four-dimensional left-sided breast cancer patient datasets, provided by the Princess Margaret Cancer Centre in Toronto, Canada. Model implementation was done using a combination of C++ and CPLEX. Instances of \eqref{RDAOMod} were run on a cluster with 125 GB RAM, while the considerably smaller CPG heuristic was run on a 4--16 GB RAM allocation.

The \eqref{RDAOMod} model and the CPG heuristic are both capable of independently generating deliverable, robust IMRT treatment plans. In the presence of the warm start generation algorithm (Algorithm \ref{warmStA}), they can be combined into the warm-started \eqref{RDAOMod}, which provides a third plan generation option. 
The remainder of this section contains a performance analysis of each of the three methods of plan generation and is organized as follows. 
Section~\ref{sec:dataSpec} outlines the data and parameter specifications, Section~\ref{Sec:Bounds} shows the CPG heuristic in action, Section~\ref{Sec:Compute} reports on the outcomes of running the \eqref{RDAOMod} model with no warm start. Finally, Section~\ref{Sec:Compute2} dives into the results of running the \eqref{RDAOMod} model with the CPG heuristic used as a warm start.

\subsection{Data Specification and Parameterization}
\label{sec:dataSpec}
Breast cancer patients A, B, C, D and E were each prescribed a dose of 42.4 Gy, to be delivered by two-beam (angle) tangential IMRT treatment. The goal of the treatment is to irradiate the target $\mathcal{T}$, which in this case is the whole-breast volume, while avoiding any unnecessary dose to the healthy organs $\mathcal{H}$, or the exposed region of the heart. 
Each patients' beam setup and treatment region features are summarized in Table~\ref{tbl:patientSum}.

\begin{table}[htb]
  \centering
  \small
  \setlength{\tabcolsep}{6pt}
  \begin{tabular}{@{\extracolsep{5pt}}ccccccc@{}}%{ccccc}
%\specialrule{.1em}{.05em}{.05em}
\toprule
& \multicolumn{2}{c}{\textbf{Beam}} & \multicolumn{2}{c}{\textbf{Target Voxels}} & \multicolumn{2}{c}{\textbf{Heart Voxels}} \\
    \cmidrule{2-3}\cmidrule{4-5}\cmidrule{6-7}
    %  \cline{2-7}
\textbf{} & Resolution & Beamlets $|\mathcal{B}|$  & Original & Sampled &  Original & Reduced \\
  %\cline{2-7} 
 %\specialrule{.1em}{.05em}{.05em} %
 \midrule%\hline
\bf Patient A &  $46 \times 26$ & 2,392 & 73,441 & 2,296 & 53,050 & 8,930  \\ %\hline
\bf Patient B & $40 \times 19$ & 1,520 &  33,592 & 1,050 & 42,342 & 2,898  \\
\bf Patient C & $ 46 \times 23 $ & 2,116 & 101,354   & 3,168 & 44,210  & 5,435  \\
\bf Patient D & $ 44 \times 22 $ & 1,936 &  56,923 & 1,779 & 51,709  & 2,465  \\
\bf Patient E & $ 36 \times 25$ & 1,800 &  70,050  & 2,190 & 52,847 & 3,812  \\
\bottomrule
%\specialrule{.1em}{.05em}{.05em}   
\end{tabular}
  \caption{Patient cardinality information for the five tangential breast cancer studies.}
  \label{tbl:patientSum} 
\end{table}

Before running the models, the decision spaces were pruned, first by omitting heart voxels that received zero or negligible dose from the beam, then by sampling every $32^{\text{nd}}$ target voxel; a method previously shown to maintain plan quality and compare to the granularity achieved by typical planning systems \citep{chan2014robust}. 
The final number of voxels are listed in Table~\ref{tbl:patientSum} and the datasets are included in this paper's online supplement.

The model's objective parameters for the target and heart, $c_\mathcal{T}$ and $c_\mathcal{H}$, were set to $0.7$ and $0.3$, respectively, as these values were empirically found to produce a good balance between the clinical objectives of conformity and heart sparing. 
Similarly, to compromise between plan quality and aperture flexibility, $|\mathcal{A}| = 6$ 
was selected as the fixed number of apertures for all generated plans.

In radiation therapy for breast cancer, when the patient's lungs expand and contract throughout each regular breathing cycle, deformation occurs in both the heart and the breast. To quantify the extent of breathing motion, patients may undergo a four-dimensional (4D) computed tomography (CT) scan, capturing the state of organs at various phases in the breathing cycle. These 4D-CT images are comprised of sets of 3D-CT images, wherein each image corresponds to a specific patient breathing phase.

For each patient, the heart and target in every 3D-CT image within the 4D scan were delineated separately, using RayStation (version 3.99.0.8, RaySearch Laboratories AB, Stockholm, Sweden).
Voxels were then tracked across each of the scans using the hybrid deformable image registration technique discussed in \cite{weistrand2015anaconda}, \cite{kadoya2016multi} and \cite{dir2017whitepaper}, which creates a one-to-one mapping for each voxel at all phases. The dose-influence matrix at each phase $i$, $D_{v,b,i}$, was then calculated using the research version of RayStation.
To account for breathing motion uncertainty, a set of upwards and downwards deviations from the nominal (expected) breathing proportions
were used to form the polyhedral uncertainty set \eqref{PDef}.
For each patient, $|\mathcal{I}| = 5$ breathing phases, sequenced from full inhale to full exhale, were considered. 
Nominal breathing proportions were set to $\textbf{p} = \lbrack 0.125, 0.125, 0.125, 0.125, 0.5\rbrack$,
with associated uncertainty set $\underline{\textbf{p}}=\bar{\textbf{p}}=\textbf{0.1}$ (i.e., the range around each nominal value is $\pm$0.1) used for robust planning,
based on the literature \citep{mahmoudzadeh2015robust}. Nominal plans are recovered using the singular $\mathcal{P}$ obtained when $\underline{\textbf{p}}=\bar{\textbf{p}}=\textbf{0}$. Unless indicated otherwise, the upcoming sections uses the nominal \textbf{p} value to evaluate model outputs.

In order to evaluate both clinical and run-time impacts of all model variants, the focus was placed on three key, model features; robustness, continuity and deliverability. As a result, models were run with nominal and robust uncertainty sets;
with and without the vertical and horizontal continuity constraints (which are activated/deactivated together); as well as with and without any deliverability constraints (to provide a baseline) for a total of six model variants: 1) FMO, 2) RFMO, 3) DAO, 4) DAO-C, 5) RDAO, and 6) RDAO-C. The ``-C'' indicates variants are run with continuity constraints. 

\begin{table}[h]
\centering 
\small
  \begin{tabular}{lcccc}
  \toprule
  %\specialrule{.1em}{.05em}{.05em}%\hline
   Model & Base Model & 
   %$|\mathcal{P}| > 1$ 
   Uncertainty
   & Deliverability & Continuity \\ \specialrule{.1em}{.3em}{.3em}%\midrule%\hline
   FMO & \checkmark & \xmark & \xmark & \xmark   \\
   RFMO & \checkmark & \checkmark & \xmark  & \xmark  \\
   DAO & \checkmark & \xmark & \checkmark & \xmark  \\
   RDAO & \checkmark & \checkmark &  \checkmark & \xmark  \\
   DAO-C & \checkmark & \xmark & \checkmark & \checkmark  \\
   RDAO-C & \checkmark & \checkmark & \checkmark & \checkmark \\
   \bottomrule
     %\specialrule{.1em}{.05em}{.05em}%\hline
   \end{tabular}
  \caption{Constraints sets included in each model variant.}
  \label{tbl:modelNames}
\end{table}

A breakdown of the constraints used in each model variant is found in Table \ref{tbl:modelNames}, with the complete mathematical overview presented in 
%Online 
Appendix 
%B.
\ref{Sec:fullMod}. 
Note that the Section~\ref{Sec:Aperture} symmetry-breaking constraints are included in the deliverability constraints.

Table~\ref{tbl:probSize} contains size information for all six model variants.
Recall that the FMO acts as the baseline models for DAO and DAO-C, while RFMO does the same for its mixed-integer counterparts.
It is evident that the largest jump in both constraints and variables comes from adding DAO requirements. The robust counterpart adds a relatively large number of additional variables and constraints, but none of the added variables are binary, while the continuity adds an even larger number of constraints and a relatively low number of variables, however, all additional variables are binary.

\begin{table}[htb]
  \centering
  \small
    \setlength{\tabcolsep}{6.5pt}
  \begin{tabular}{clcccccc}
 \toprule%\specialrule{.1em}{.05em}{.05em}% \hline
   \multicolumn{2}{c}{} & \multicolumn{6}{c}{\textbf{Model}} \\      \cmidrule{3-8}
 \multicolumn{1}{c}{\textbf{}} & \multicolumn{1}{l}{\textbf{Feature}} & FMO & DAO & DAO-C & RFMO & RDAO & RDAO-C \\
  \specialrule{.1em}{.3em}{.3em}%\hline
  \multirow{3}{*}{\textbf{Patient A}}  & Constraints & 2,296 & 87,332 & 118,148 & 13,776 & 98,812 & 129,628 \\
  & Variables & 2,392 & 57,426 & 59,082 & 16,168 & 71,202 & 72,858 \\
  & Binaries & 0 & 43,068 & 44,724 & 0 & 43,068 & 44,724 \\
     \midrule%%\specialrule{.1em}{.05em}{.05em}
   \multirow{3}{*}{\textbf{Patient B}} & Constraints & 1,050 & 54,838 & 74,998 & 6,300 & 60,088 & 80,248 \\
   & Variables & 1,520 & 36,498 & 37,938 & 7,820 & 42,798 & 44,238 \\
   & Binaries & 0 & 27,372 & 28,812 & 0 & 27,372 & 28,812 \\
    \midrule% \specialrule{.1em}{.05em}{.05em}%\hline
        \multirow{3}{*}{\textbf{Patient C}} & Constraints & 3,168 & 78,268 & 105,844 & 19,008 & 94,108  & 121,684 \\
   & Variables & 2,116 & 50,802 & 52,458 & 21,124 & 69,810 & 71,466  \\
   & Binaries & 0 & 38,100 & 39,756 & 0 & 38,100 & 39,756  \\
        \midrule%\specialrule{.1em}{.05em}{.05em}%\hline
           \multirow{3}{*}{\textbf{Patient D}} & Constraints & 1,779 & 70,447 & 95,767 & 10,674 & 79,342 & 104,662 \\
   & Variables & 1,936 & 46,482 & 48,066 & 12,610 & 57,156 & 58,740 \\
   & Binaries & 0 & 34,860 & 36,444 & 0 & 34,860 & 36,444 \\
       \midrule% \specialrule{.1em}{.05em}{.05em}%\hline
           \multirow{3}{*}{\textbf{Patient E}} & Constraints & 2,190 & 66,154 & 89,290 & 13,140 & 77,104  & 100,240 \\
   & Variables & 1,800 & 43,218 & 44,514 & 14,940 & 56,358 & 57,654 \\
   & Binaries & 0 & 32,412 & 33,708 & 0 & 32,412 & 33,708 \\
     \bottomrule%\specialrule{.1em}{.05em}{.05em}%\hline
   \end{tabular}
  \caption{Problem size of the six model variants for each patient. }
  \label{tbl:probSize}
\end{table}

\subsection{CPG Heuristic Performance}
\label{Sec:Bounds}

Running the CPG heuristic on its own yields a clinically viable IMRT plan. To gain insight into the heuristic, the output can be visualized at each step using a fluence map, which shows dose intensity from the beam of radiation's perspective. Within the map, lighter beamlets depict higher intensities and blackened beamlets have zero intensity.

\begin{figure}[h]
    \begin{subfigure}[t]{0.08\textwidth}
        \centering
         \caption*{\small \textbf{Step 1}}
        	\includegraphics[width=.98\textwidth, height = 5cm]
        	{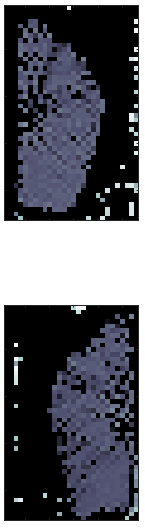}%{BEV_LogFMO_32_Pat1_2.png}          
         \label{fig:ws_s1_pat1}
    \end{subfigure} %
        \begin{subfigure}[t]{0.03\textwidth}
        \centering
	\raisebox{-7\height}{\includegraphics[width=.9\textwidth, height = .5cm]{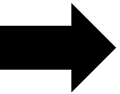}}    
	\end{subfigure} %
    \begin{subfigure}[t]{0.29\textwidth}
        \centering
              \caption*{\small \textbf{Step 2}}
        \includegraphics[width=.9\textwidth,height = 5cm]{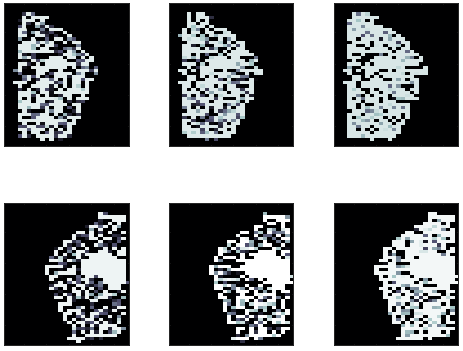}%{BEV_FMO_32_Pat1_Warm2_2.png}
         \label{fig:ws_s2_pat1}
    \end{subfigure} %
            \begin{subfigure}[t]{0.03\textwidth}
        \centering
        %\raisebox{-0.5\height}
	\raisebox{-7\height}{\includegraphics[width=.9\textwidth,height = .5cm]{arrow.png}}            
~~~~~~ \end{subfigure} %
    \begin{subfigure}[t]{.56\textwidth}
        \centering
                \caption*{\small \textbf{Step 3}}
        \includegraphics[width=.98\textwidth, height = 5cm]{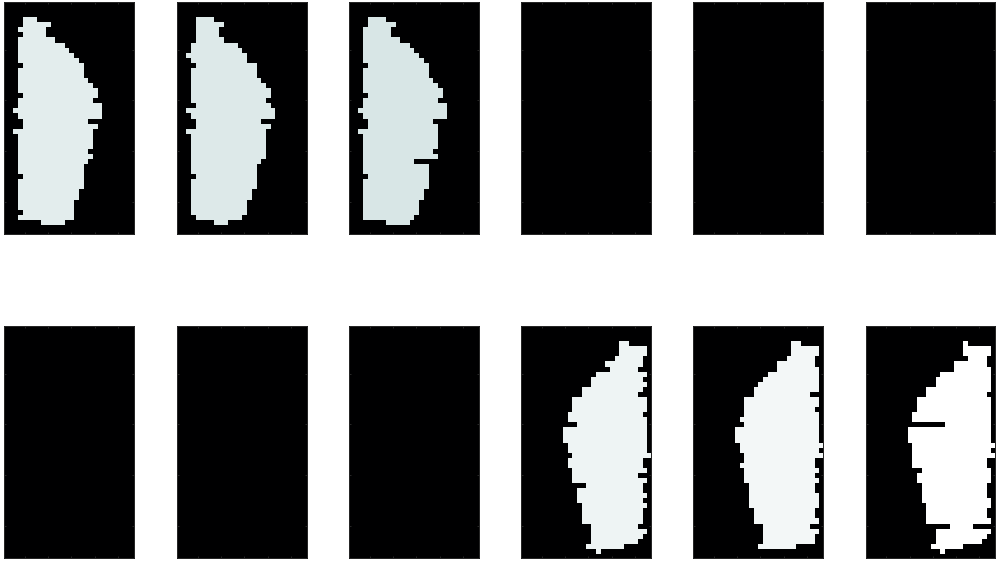}
        %BEV_FMO_32_Pat1_Warm3_CB2.png}
         \label{fig:ws_s3_pat1}
    \end{subfigure} %
    \caption{The three-step CPG algorithm output demonstrated on Patient B. A log$_{10}$ scale is applied to the beamlet intensities in Step~1, whereas Steps 2 and 3 have a linear scale.}
    \label{fig:WarmStartDemo}
\end{figure}

Figure~\ref{fig:WarmStartDemo} depicts the fluence map of the output of all three steps of the CPG heuristic. Note the model output is visually similar across all patients and model variants, so the selection of Patient B's DAO output was arbitrary. The Step~1 output is exactly the FMO intensities, since the model was DAO. Had a robust model variant (RDAO or RDAO-C) been chosen, the visually-similar RFMO output would have been presented instead. 
A log$_{10}$ scale was applied to the FMO fluence map 
because the contrast between high and low intensities
was so large that most active beamlets did not show up on the map without a change in scaling.
In the Figure~\ref{fig:WarmStartDemo} example, the highest intensity was 104,051 units, which could take on the order of weeks to deliver. In contrast, the maximum intensity seen in the following steps is a much more reasonable 56 units. 

Figure~\ref{fig:WarmStartDemo}'s Step~2 fluence map depicts a more-uniform, but still undeliverable \eqref{WSMod} plan, with only $\frac{1}{|\mathcal{A}|} = \frac{1}{2}$ of the required apertures.
Finally, Step~3 reveals a properly-sorted, deliverable six-aperture plan. 
The final plan, while not optimal,  
could be delivered without adjustment, if a clinician felt it was of sufficiently high-quality. One way to assess the quality of a CPG plan is through the comparison of its objective function value ($z^{CPG}$) to that of its Step 1 linear \eqref{RobMod} counterpart, which acts as a lower bound.

\begin{table}[htb]
  \centering
  \small
   \setlength{\tabcolsep}{9pt}
  \begin{tabular}{clccc >{\color{black}}c}
 %\specialrule{.1em}{.05em}{.05em}%\hline
 \toprule
 \multicolumn{1}{c}{} & \multicolumn{1}{c}{\textbf{Model}} & \textbf{Lower Bound} 
 & $\mathbf{z^{CPG}}$ & \textbf{\% Difference} &  $\alpha$ \\
  \specialrule{.1em}{.3em}{.3em}%\cmidrule{1-6}%
  \multirow{2}{*}{\textbf{Patient A}}  & DAO/DAO-C & 36.38  & 55.33 &  34.24 & 0.4  \\
   &RDAO/RDAO-C  & 36.89  & 56.76 & 35.00 & 0.8  \\
     \cmidrule{1-6}%\hline
   \multirow{2}{*}{\textbf{Patient B}} & DAO/DAO-C & 31.27  & 37.65 & 16.94 & 0.4 \\
   & RDAO/RDAO-C   & 31.42  & 38.22 & 17.78 & 0.4 \\
        \cmidrule{1-6}%\hline
   \multirow{2}{*}{\textbf{Patient C}} & DAO/DAO-C & 34.32 & 51.54 & 33.42  & 0.8 \\
   & RDAO/RDAO-C   & 34.51 & 51.81 &  33.40  & 0.8 \\
     \cmidrule{1-6}%\hline
   \multirow{2}{*}{\textbf{Patient D}} & DAO/DAO-C & 32.62 & 50.55 & 35.47 & 0.2 \\
   & RDAO/RDAO-C   &  33.18 & 51.39 & 35.44 & 0.2\\
     \cmidrule{1-6}%\hline
   \multirow{2}{*}{\textbf{Patient E}} & DAO/DAO-C & 33.25 & 51.72 & 35.70 & 0.4\\
   & RDAO/RDAO-C & 33.47 & 52.98 & 36.82 & 0.8 \\
     %\specialrule{.1em}{.05em}{.05em}%
     \bottomrule
   \end{tabular}
  \caption{Table of objective function values ($z^*$) of the linear FMO/RFMO models (lower bounds), compared to the CPG objective function values, and the resulting worst-case optimality gap.}
  \label{tbl:bounds}
\end{table}

Table~\ref{tbl:bounds} contains the CPG heuristic and linear lower bounds on the objective function value for each variant of the model. See 
%Online 
Appendix 
\ref{app:alphaTune} 
%C 
for the tuning of parameter $\alpha$.
Note that the DAO and DAO-C, as well as RDAO and RDAO-C, have been merged
in the table because none of the DAO or RDAO plans generated had discontinuities, likely due to the geometry of the target at hand.
That means that the plans are identical to the those produced for DAO-C and RDAO-C, respectively.

The second-to-last column in Table~\ref{tbl:bounds} shows the optimality gap between the lower bound obtained from \eqref{RobMod} and the CPG heuristic plan. 
This value, which is on average 34.6\%, 17.4\%, 33.4\%, 35.5\% and 36.3\% for patients A to E, provides an upper bound on the true optimality gap of the CPG plans. It likely overestimates the true gap, since the linear models relax deliverability constraints and are therefore not feasible plans. For the RDAO models, the plan also depends on the size of the uncertainty set. See 
%Online 
Appendix 
\ref{app:uncertain}
%D 
for a brief analysis of this relationship.

\subsection{Basic RDAO Model Performance}
\label{Sec:Compute}
To test the exact holistic model, all four MIP variants of \eqref{RDAOMod} were run with no warm start, then compared to CPG heuristic output for each case study. Since \eqref{RDAOMod} is a very large-scale problem, each model variant was allocated a total runtime of three days. Table~\ref{tbl:AllDAOResults} shows each model's objective function value ($z$), the time at which it was obtained and its corresponding CPLEX-reported optimality gap at the first feasible integer solution (\textit{First Incumbent}) and after the three-day runtime (\textit{Best Incumbent}). 
Finally, the CPG objective function values, runtimes and worst-case optimality gaps, are provided. 

\begin{table}[]
\centering
\small
\begin{tabular}{@{}clccccc>{\color{black}}cccc@{}}
\toprule
 &  & \multicolumn{3}{c}{\textbf{First Incumbent}} & \multicolumn{3}{c}{\textbf{Best Incumbent}} & \multicolumn{3}{c}{\textbf{CPG Heuristic}} \\ \cmidrule(l){3-5} \cmidrule(l){6-8} \cmidrule(l){9-11}
                                    & \textbf{Model} & $z^{Inc}$ & Gap    & Time          & $z^{Best}$ & Gap    & Time  & $z^{CPG}$     & Gap    & Time          \\ \specialrule{.1em}{.3em}{.3em}%\bottomrule
\multirow{4}{*}{\textbf{Patient A}} & DAO            & 59.1      & 38.1\% & \textbf{0.21} & 58.32      & 37.3\% & 9.56  & \textbf{55.3} & 34.2\% & 0.33          \\
                                    & DAO-C          & -         & -      & -             & -          & -      & -     & \textbf{55.3} & 34.2\% & \textbf{0.33} \\
                                    & RDAO           & 60.5      & 38.6\% & 2.85          & 60.46      & 38.6\% & 2.85  & \textbf{56.8} & 35.0\% & \textbf{0.81} \\
                                    & RDAO-C         & -         & -      & -             & -          & -      & -     & \textbf{56.8} & 35.0\% & \textbf{0.81} \\ \midrule
\multirow{4}{*}{\textbf{Patient B}} & DAO            & 57.4      & 45.5\% & \textbf{0.04} & 56.88      & 45.0\% & 1.34  & \textbf{37.7} & 16.9\% & 0.05          \\
                                    & DAO-C          & 3839.2    & 99.2\% & 5.02          & 3046.26    & 99.0\% & 12.64 & \textbf{37.7} & 16.9\% & \textbf{0.05} \\
                                    & RDAO           & 57.6      & 45.4\% & 0.21          & 57.58      & 45.4\% & 0.21  & \textbf{38.2} & 17.8\% & \textbf{0.14} \\
                                    & RDAO-C         & -         & -      & -             & -          & -      & -     & \textbf{38.2} & 17.8\% & \textbf{0.14} \\ \midrule
\multirow{4}{*}{\textbf{Patient C}} & DAO            & 71.0      & 51.6\% & \textbf{0.32} & 64.59      & 46.7\% & 32.35 & \textbf{51.5} & 33.4\% & 0.76          \\
                                    & DAO-C          & 79.1      & 56.5\% & 41.61         & 76.70      & 55.2\% & 41.85 & \textbf{51.5} & 33.4\% & \textbf{0.76} \\
                                    & RDAO           & 71.1      & 51.5\% & 5.07          & 71.13      & 51.5\% & 5.07  & \textbf{51.8} & 33.4\% & \textbf{1.97} \\
                                    & RDAO-C         & -         & -      & -             & -          & -      & -     & \textbf{51.8} & 33.4\% & \textbf{1.97} \\ \midrule
\multirow{4}{*}{\textbf{Patient D}} & DAO            & 64.2      & 48.5\% & \textbf{0.14} & 59.87      & 44.8\% & 12.26 & \textbf{50.6} & 35.5\% & 0.32          \\
                                    & DAO-C          & -         & -      & -             & -          & -      & -     & \textbf{50.6} & 35.5\% & \textbf{0.32} \\
                                    & RDAO           & 64.5      & 48.0\% & \textbf{1.20} & 64.47      & 48.0\% & \textbf{1.20}  & \textbf{51.3} & 35.3\% & 1.40          \\
                                    & RDAO-C         & -         & -      & -             & -          & -      & -     & \textbf{51.3} & 35.3\% & \textbf{1.40} \\ \midrule
\multirow{4}{*}{\textbf{Patient E}} & DAO            & 64.8      & 48.5\% & \textbf{0.12} & 64.79      & 48.5\% & \textbf{0.12}  & \textbf{51.7} & 35.7\% & 0.24          \\
                                    & DAO-C          & 3816.4    & 99.1\% & 14.61         & 3529.40    & 99.1\% & 63.49 & \textbf{51.7} & 35.7\% & \textbf{0.24} \\
                                    & RDAO           & 65.6      & 48.7\% & 1.35          & 65.57      & 48.7\% & 1.35  & \textbf{53.0} & 36.8\% & \textbf{0.34} \\
                                    & RDAO-C         & -         & -      & -             & -          & -      & -     & \textbf{53.0} & 36.8\% & \textbf{0.34} \\ \bottomrule
\end{tabular}
 \caption{The first and best incumbent solutions of the 4 integer variants of the \eqref{RDAOMod} model for all patients (run for 3-days), compared to the heuristic CPG plan objective. A dash indicates no solution is found, bold values show the best (lowest) optimality gaps and runtimes (given in hours, wall-clock time).}
  \label{tbl:AllDAOResults}
\end{table}

It is evident from Table~\ref{tbl:AllDAOResults} that the CPG heuristic dominates, not only the first incumbent values found by the solver for \eqref{RDAOMod}, but also every best CPLEX-obtained incumbent value across the board. 
The discrepancies are most pronounced in the cases with continuity constraints, which make the models more complex to solve.
In particular, none of the models with robust continuity constraints (RDAO-C) find a feasible incumbent solution within the three day time limit, and DAO-C proves equally difficult for patients A and D.
In contrast, the CPG heuristic finds solutions to every instance with optimality gaps ranging between 16.9\% and 36.8\%. The remaining DAO-C models hardly fared better, with the solver proposing first incumbents objective function values that are orders of magnitude greater than those produced by the CPG heuristic for patients B and D, with only patient C finding a reasonable solution within the allocated time period.

The runtime information in Table~\ref{tbl:AllDAOResults} also appears to confirm the CPG heuristic's dominance over finding a first incumbent solution in larger cases. While the CPG runtimes varied from under three minutes, to just under two hours, many of the \eqref{RDAOMod} variants spent days finding plans that were ultimately lower quality than the CPG.
Even in the six cases where the incumbent is found faster using \eqref{RDAOMod}, the CPG model catches up in a matter of minutes, significantly overtaking the solution in terms of quality.

Table \ref{tbl:PlanQual} provides information about the basic-\eqref{RDAOMod}-generated and CPG-generated plan alternatives. Within the table, all plans were normalized to provide 95\% of the desired dose to 95\% of the target volume; a practice commonly used in clinic. Given this normalization, the average and highest dose-to target ($\mathcal{T}^{ave}$ and $\mathcal{T}^{max}$) as well as the average and highest dose to the heart ($\mathcal{H}^{ave}$ and $\mathcal{H}^{max}$) are reported for each plan. 

\begin{table}[h]
\centering
\small
\begin{tabular}{@{\extracolsep{4pt}}clcccccccc@{}}
\toprule
                                    &        & \multicolumn{4}{c}{Best Basic (RDAO) Plan}      & \multicolumn{4}{c}{CPG Heuristic Plan}                       \\  \cmidrule(l){3-6} \cmidrule(l){7-10}
 &
  \textbf{Model} &
  $\mathcal{T}^{ave}$ &
  $\mathcal{T}^{max}$ &
  $\mathcal{H}^{ave}$ &
  $\mathcal{H}^{max}$ &
  ~$\mathcal{T}^{ave}$~~&
  ~$\mathcal{T}^{max}$~&
  $\mathcal{H}^{ave}$&
  $\mathcal{H}^{max}$~ \\ \midrule
\multirow{4}{*}{\textbf{Patient A}} & DAO    & 56.8   & 78.4    & 5.0          & 59.8          & \textbf{51.8} & \textbf{65.9} & \textbf{4.2} & \textbf{53.1} \\
                                    & DAO-C  & -      & -       & -            & -             & \textbf{51.8} & \textbf{65.9} & \textbf{4.2} & \textbf{53.1} \\
                                    & RDAO   & 57.8   & 79.2    & 5.3          & 60.6          & \textbf{55.5} & \textbf{74.1} & \textbf{5.0} & \textbf{57.8} \\
                                    & RDAO-C & -      & -       & -            & -             & \textbf{55.5} & \textbf{74.1} & \textbf{5.0} & \textbf{57.8} \\ \midrule
\multirow{4}{*}{\textbf{Patient B}} & DAO    & 54.4   & 72.7    & 1.3          & 51.6          & \textbf{44.4} & \textbf{48.8} & \textbf{0.6} & \textbf{42.1} \\
                                    & DAO-C  & 1835.6 & 26769.2 & 32.5         & 15193.9       & \textbf{44.4} & \textbf{48.8} & \textbf{0.6} & \textbf{42.1} \\
                                    & RDAO   & 54.6   & 72.9    & 1.4          & 52.1          & \textbf{45.0} & \textbf{49.4} & \textbf{0.5} & \textbf{42.4} \\
                                    & RDAO-C & -      & -       & -            & -             & \textbf{45.0} & \textbf{49.4} & \textbf{0.5} & \textbf{42.4} \\ \midrule
\multirow{4}{*}{\textbf{Patient C}} & DAO    & 67.5   & 105.9   & \textbf{2.2} & 61.2          & \textbf{51.5} & \textbf{65.6} & 2.3          & \textbf{54.6} \\
                                    & DAO-C  & 69.9   & 105.5   & 2.8          & 84.4          & \textbf{51.5} & \textbf{65.6} & \textbf{2.3} & \textbf{54.6} \\
                                    & RDAO   & 67.6   & 101.2   & 3.3          & 63.4          & \textbf{51.9} & \textbf{66.2} & \textbf{2.4} & \textbf{55.1} \\
                                    & RDAO-C & -      & -       & -            & -             & \textbf{51.9} & \textbf{66.2} & \textbf{2.4} & \textbf{55.1} \\ \midrule
\multirow{4}{*}{\textbf{Patient D}} & DAO    & 59.6   & 85.2    & 0.7          & 50.3          & \textbf{46.1} & \textbf{52.7} & \textbf{0.4} & \textbf{41.6} \\
                                    & DAO-C  & -      & -       & -            & -             & \textbf{46.1} & \textbf{52.7} & \textbf{0.4} & \textbf{41.6} \\
                                    & RDAO   & 59.2   & 84.2    & 0.9          & 54.3          & \textbf{48.5} & \textbf{55.7} & \textbf{0.4} & \textbf{43.5} \\
                                    & RDAO-C & -      & -       & -            & -             & \textbf{48.5} & \textbf{55.7} & \textbf{0.4} & \textbf{43.5} \\ \midrule
\multirow{4}{*}{\textbf{Patient E}} & DAO    & 65.0   & 93.4    & 1.6          & 56.4          & \textbf{51.3} & \textbf{64.8} & \textbf{1.1} & \textbf{54.6} \\
                                    & DAO-C  & 1373.0 & 33006.1 & 9.2          & 2888.9        & \textbf{51.3} & \textbf{64.8} & \textbf{1.1} & \textbf{54.6} \\
                                    & RDAO   & 65.7   & 94.4    & 1.6          & \textbf{57.1} & \textbf{53.9} & \textbf{70.9} & \textbf{1.3} & 59.2          \\
                                    & RDAO-C & -      & -       & -            & -             & \textbf{53.9} & \textbf{70.9} & \textbf{1.3} & \textbf{59.2} \\ \bottomrule
\end{tabular}
 \caption{Normalized plans from the integer \eqref{RDAOMod} variants using a basic solver (for 3-days) compared to the heuristic CPG values. Bolded values are the best (lowest) average or maximum doses to the target $\mathcal{T}$, or heart $\mathcal{H}$.}
  \label{tbl:PlanQual}
\end{table}
 
Based on these quality assessment metrics, it is evident that the CPG heuristic continues to outperform 3-days of solver runtime. With the exception of Patient C's DAO average heart dose and Patient E's RDAO maximum heart dose, all values output were lower for the CPG plan. It should be noted that when normalization was applied to the robust cases, 95\% of the minimum dose in the generated plan was used, rather than 95\% of the prescribed dose, in an effort to maintain the robust effect. There could, however, still potentially be some erosion of the robust guarantee due to interpolation in these instances.

\begin{figure} [htb]
    \centering
    \begin{subfigure}[b]{0.49\textwidth}
        \centering
        \includegraphics[width=.9\textwidth]{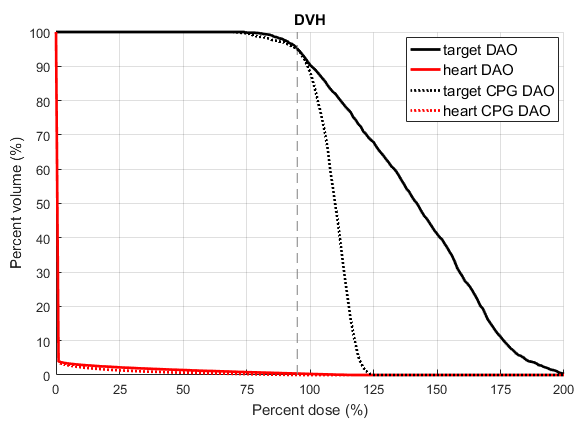}
        \caption{DVH of DAO plan}
        \label{fig:RDAODVH}
    \end{subfigure}%
    ~~~
    \begin{subfigure}[b]{0.49\textwidth}
        \centering
        \includegraphics[width=.9\textwidth]{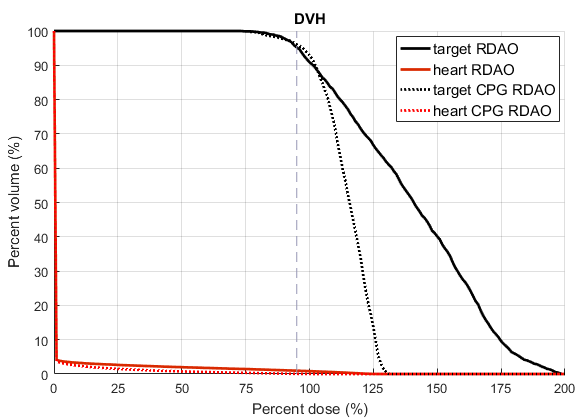}
        \caption{DVH of RDAO plan}
        \label{fig:DAODVH}
    \end{subfigure}
    \caption{The DVH of Patient D's DAO and RDAO plans generated using CPLEX (solid line), compared to the DAO and RDAO plan generated using the CPG heuristic (dotted line).}
     \label{fig:Pat3DAODVH}
\end{figure}

In clinic, the quality of a plan is often assessed using a tool called a dose-volume histogram (DVH), which depicts the percentage of the prescribed dose received by a particular volume of an organ in the region of interest. Using this tool, the plans above may be depicted and compared. For example, Patient D's normalized DAO and RDAO plans are shown in Figure \ref{fig:Pat3DAODVH}, obtained by running the \eqref{RDAOMod} models (solid line) and the CPG heuristic (dotted line). Note that the trend of quality gain with the CPG heuristic continues to be seen, as a much more conformal plan appears to be generated, in agreement with Table \ref{tbl:PlanQual}, above. The other patients (omitted, due to being visually similar to Patient D) follow this same trend.

\subsection{Warm-Started RDAO Model Performance}
\label{Sec:Compute2}
In this section, the quality improvement of the plans generated with the CPG heuristic warm starting \eqref{RDAOMod}, over the solutions of \eqref{RDAOMod} alone is demonstrated, followed by a discussion of the warm-started \eqref{RDAOMod} results.

Table~\ref{tbl:WarmStartResults} presents the objective function values and optimality gaps of the best solutions found using the \eqref{RDAOMod} variants both alone and with the CPG algorithm as a warm start. The table is partitioned horizontally, based on whether or not the model was run with continuity constraints (-C), seeing as CPLEX struggled to find reasonable incumbents in these cases, making the difference far more pronounced (with the exception of Patient C's DAO-C).
Note that the objective function values of the non-warm-started (No WS) \eqref{RDAOMod} model after three days, $z^{Best}$ No WS, were
previously reported in Table~\ref{tbl:AllDAOResults}, under \textit{Best Incumbent}, and are duplicated here for reference.
The warm-started (WS) models consistently outperform the No WS models. As such, the percent reduction in the optimality gap of the warm-started \eqref{RDAOMod}, compared to the non-warm started best gap is provided in the final row of each section of the Table~\ref{tbl:WarmStartResults}. These values range from a low-end improvement of 11.6\% to a high-end of 75.8\%, with no continuity constraints, and from 43.8\% to 84.6\% with continuity constraints. 
\begin{table}[]
\centering
\resizebox{\textwidth}{!}{%
\begin{tabular}{@{}llcccccccccc@{}}
\toprule
 &  & \multicolumn{2}{c}{\textbf{Patient A}} & \multicolumn{2}{c}{\textbf{Patient B}} & \multicolumn{2}{c}{\textbf{Patient C}} & \multicolumn{2}{c}{\textbf{Patient D}} & \multicolumn{2}{c}{\textbf{Patient E}} \\ %\cmidrule(l){3-12} 
\multirow{6}{*}{\rotatebox{90}{\textbf{No Continuity}}} & \textbf{} & DAO & RDAO & DAO & RDAO & DAO & RDAO & DAO & RDAO & DAO & RDAO \\  \cmidrule(l){2-2} \cmidrule(l){3-4} \cmidrule(l){5-6}\cmidrule(l){7-8}\cmidrule(l){9-10}\cmidrule(l){11-12} %\midrule
 & $z^{Best}$ No WS & 58.3 & 60.5 & 56.9 & 57.6 & 64.6 & 71.1 & 59.9 & 64.5 & 64.8 & 65.6 \\
 & Gap (\%) & 37.3 & 38.6 & 45.0 & 45.4 & 46.7 & 51.5 & 44.8 & 48.0 & 48.5 & 48.7 \\ \cmidrule(l){2-12} 
 & $z^{Best}$ CPG WS & 54.5 & 56.6 & 35.2 & 37.6 & 50.4 & 51.8 & 50.6 & 51.3 & 51.7 & 53.0 \\
 & Gap (\%) & 32.8 & 34.1 & 10.9 & 16.2 & 31.6 & 33.2 & 30.7 & 33.4 & 29.2 & 36.2 \\ \cmidrule(l){2-12} 
 & \textbf{Reduction} & \textbf{12.1\%} & \textbf{11.6\%} & \textbf{75.8\%} & \textbf{64.3\%} & \textbf{32.3\%} & \textbf{35.5\%} & \textbf{31.5\%} & \textbf{30.5\%} & \textbf{39.8\%} & \textbf{25.7\%} \\ \toprule
  &  & \multicolumn{2}{c}{\textbf{Patient A}} & \multicolumn{2}{c}{\textbf{Patient B}} & \multicolumn{2}{c}{\textbf{Patient C}} & \multicolumn{2}{c}{\textbf{Patient D}} & \multicolumn{2}{c}{\textbf{Patient E}} \\
\multirow{6}{*}{\rotatebox{90}{\textbf{-C Constraints}}} & \textbf{} & DAO-C & RDAO-C & DAO-C & RDAO-C & DAO-C & RDAO-C & DAO-C & RDAO-C & DAO-C & RDAO-C \\   \cmidrule(l){2-2} \cmidrule(l){3-4} \cmidrule(l){5-6}\cmidrule(l){7-8}\cmidrule(l){9-10}\cmidrule(l){11-12}  %\cmidrule(l){2-12} 
 & $z^{Best}$ No WS & - & - & 3046.3 & - & 76.7 & - & - & - & 3529.4 & - \\
 & Gap (\%) & 100.0 & 100.0 & 99.0 & 100.0 & 55.2 & 100.0 & 100.0 & 100.0 & 99.1 & 100.0 \\ \cmidrule(l){3-12} 
 & $z^{Best}$ CPG WS & 54.7 & 56.6 & 37.0 & 37.8 & 49.9 & 51.8 & 49.2 & 51.3 & 49.5 & 52.9 \\
 & Gap (\%) & 33.1 & 34.2 & 15.3 & 16.8 & 31.0 & 33.2 & 31.9 & 33.4 & 32.5 & 36.2 \\ \cmidrule(l){2-12} 
 & \textbf{Reduction} & \textbf{66.9\%} & \textbf{65.8\%} & \textbf{84.6\%} & \textbf{83.2\%} & \textbf{43.8\%} & \textbf{66.8\%} & \textbf{68.1\%} & \textbf{66.6\%} & \textbf{67.2\% } & \textbf{63.8\%} \\ \bottomrule
\end{tabular}%
}
\caption{Final objective function values and optimality gaps for each patient and model variant, compared when run without and with the CPG warm start.}
\label{tbl:WarmStartResults}
\end{table}

\begin{figure}[htbp]
    \centering
    \includegraphics[%trim={0.5cm 0.5cm 1cm 0.5cm},clip,
    width=\textwidth]{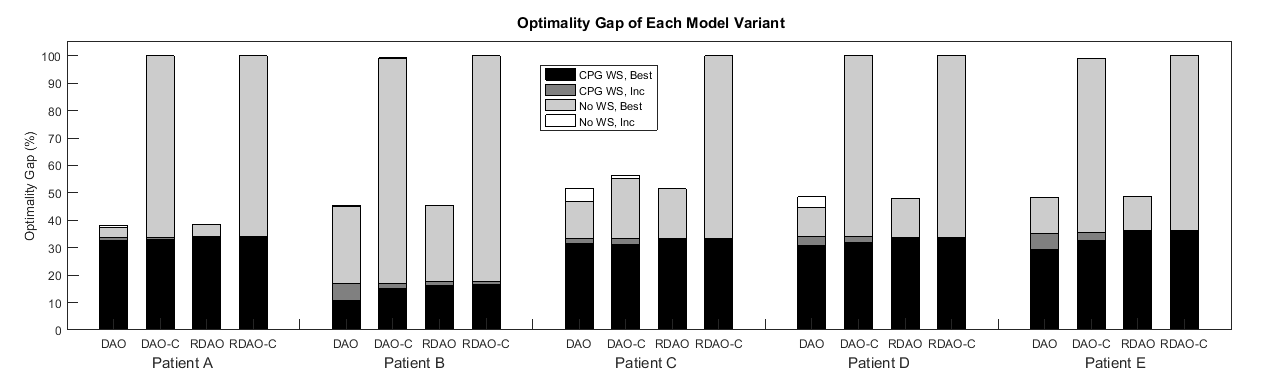}
    \caption{The cplex-reported optimality gap present at the first incumbent and 3-day stages, for both warm started and independently-run RDAO variants.}
    \label{fig:optGaps}
\end{figure}

Figure~\ref{fig:optGaps} graphically depicts the variation in optimality gaps from the first and best incumbent plans reported in Tables~\ref{tbl:AllDAOResults} and \ref{tbl:WarmStartResults}. It can be observed that the gaps may be sorted in a monotonically increasing order, with the warm-started \eqref{RDAOMod} after three days (\textit{CPG WS, Best}) being the smallest, or closest to optimality, followed by its first incumbent, or CPG heuristic output (\textit{CPG WS, Inc}). The next closest to optimal is the best independent \eqref{RDAOMod} after seven days (\textit{No WS, Best}), followed by (and sometimes tied with) the first incumbent solution found by the solver (\textit{No WS, Inc}). 
Note that the difference between the non-warm started first and best incumbent solutions are often quite small (in the 13/20 cases where an incumbent is found at all). We hypothesize that with the better starting solution, $z^{CPG}$, the MIP solver is able to prune more branches, allowing it to more effectively traverse the solution space.

\begin{figure}[htbp]
% Make the title box
     \begin{subfigure}[t]{0.05\textwidth}
          \hspace{0.05\textwidth}
     \end{subfigure}% 	
      \begin{subfigure}[t]{0.4\textwidth}
           \centering
     	  \textbf{Full Fluence Map}
     \end{subfigure} %
%          \begin{subfigure}[t]{0.05\textwidth}
%          \hspace{0.05\textwidth}
%     \end{subfigure}% 	 
     \begin{subfigure}[t]{0.55\textwidth}
        \centering
	\textbf{Aperture Breakdown}
     \end{subfigure}
    	\par\medskip

       \begin{subfigure}[t]{0.05\textwidth}
              \flushright
     	\raisebox{1.35\height}{\rotatebox{90}{\textbf{DAO}}}
    \end{subfigure}% 	
    \begin{subfigure}[t]{0.4\textwidth}
        \centering
        \includegraphics[width=.75\textwidth, height=4cm]{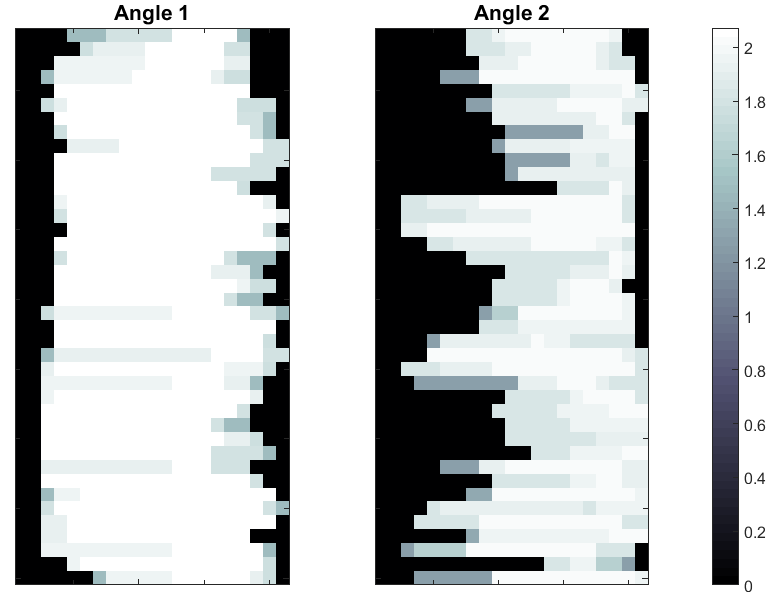}
%        \caption{}%Combined FM of NC DAO}
         \label{fig:FMDAONC}
    \end{subfigure} %
%         \begin{subfigure}[t]{0.05\textwidth}
%          \hspace{0.05\textwidth}
%     \end{subfigure}% 	 	
     \begin{subfigure}[t]{0.55\textwidth}
        \centering
        \includegraphics[width=.85\textwidth, height=4cm]{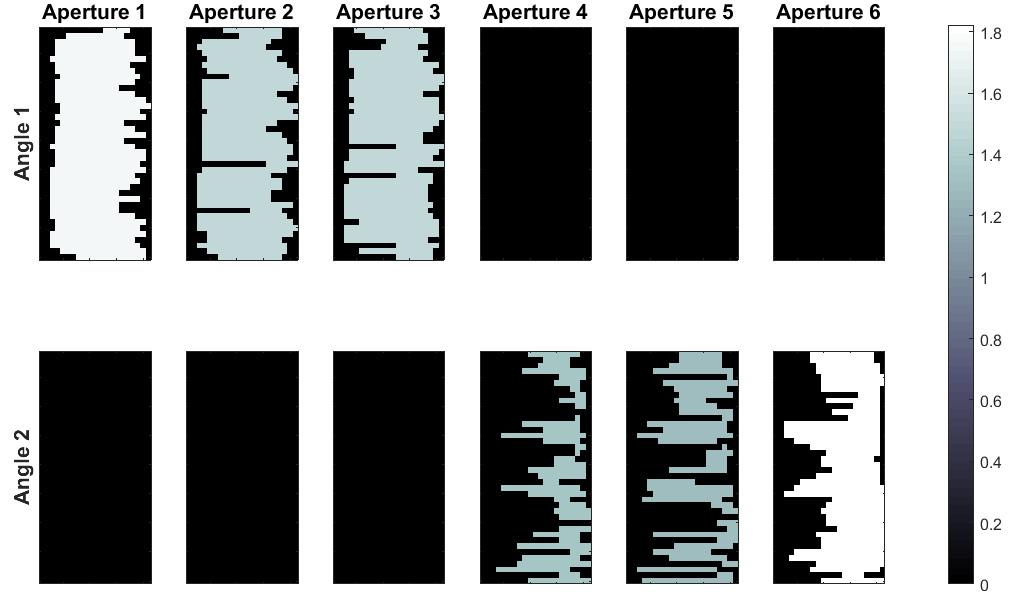}
%        \caption{}%FM of the DAO (no continuity constraints)}
         \label{fig:BEVDAONC}
    \end{subfigure}%    

    	\par\bigskip
	
       \begin{subfigure}[t]{0.05\textwidth}
       \flushright
     \raisebox{.75\height}{\rotatebox{90}{\textbf{DAO-C}}} 
    \end{subfigure}% 	
        \begin{subfigure}[t]{0.4\textwidth}
        \centering
        \includegraphics[width=.75\textwidth, height=4cm]{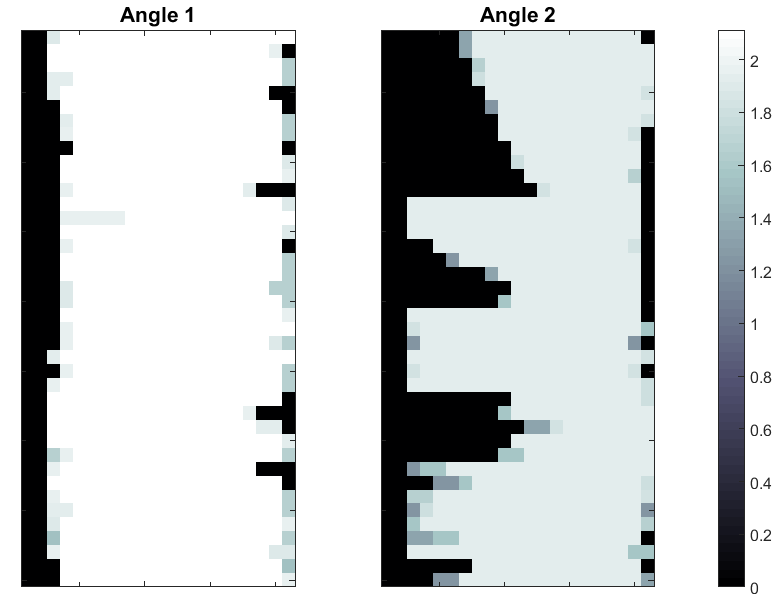}
%        \caption{}%Combined FM of NC DAO}
         \label{fig:FMDAO}
    \end{subfigure} %
%     \begin{subfigure}[t]{0.05\textwidth}
%          \hspace{0.05\textwidth}
%     \end{subfigure}% 	 	
            \begin{subfigure}[t]{0.55\textwidth}
        \centering
        \includegraphics[width=.85\textwidth, height=4cm]{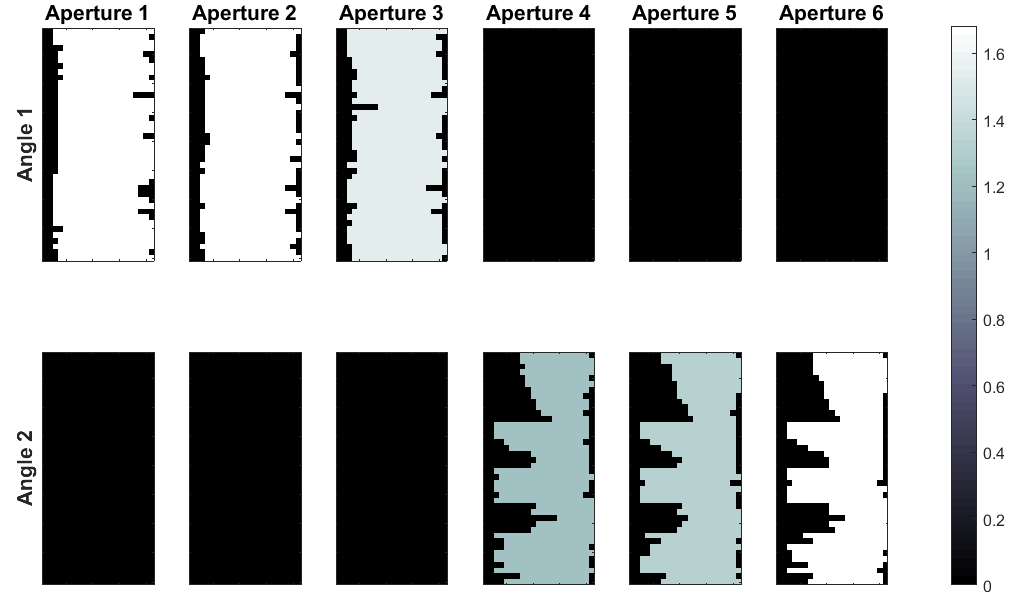}
%        \caption{}%FM of the DAO (no continuity constraints)}
         \label{fig:BEVDAO}
    \end{subfigure}
    	\par\bigskip
	
	\begin{subfigure}[t]{0.05\textwidth}
	       \flushright
      \raisebox{.9\height}{\rotatebox{90}{\textbf{RDAO}}}
    \end{subfigure}% 	
    \begin{subfigure}[t]{0.4\textwidth}
        \centering
        \includegraphics[width=.75\textwidth, height=4cm]{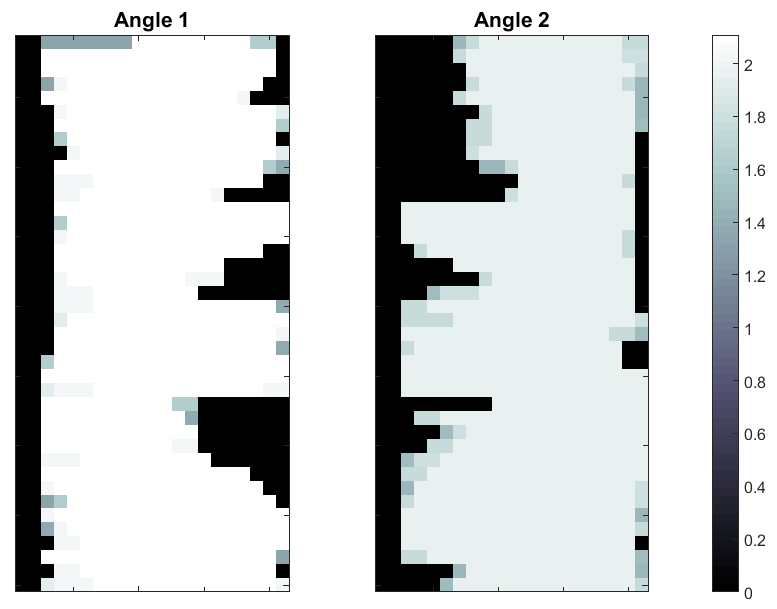}
%        \caption{}%Combined FM of NC DAO}
         \label{fig:FMRDAONC}
    \end{subfigure} %
%         \begin{subfigure}[t]{0.05\textwidth}
%          \hspace{0.05\textwidth}
%     \end{subfigure}% 	 	
     \begin{subfigure}[t]{0.55\textwidth}
        \centering
        \includegraphics[width=.85\textwidth, height=4cm]{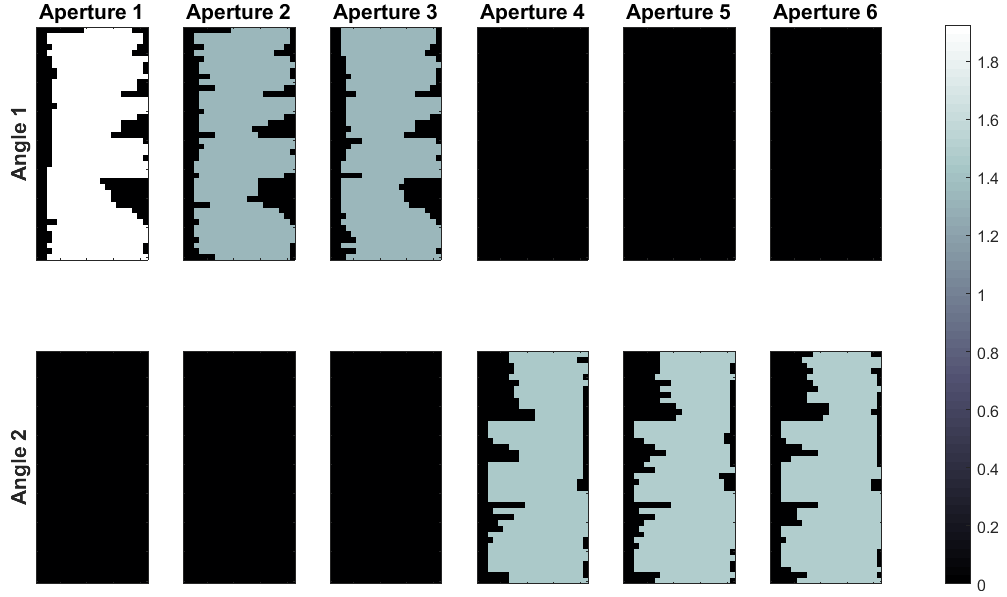}
%        \caption{}%FM of the DAO (no continuity constraints)}
         \label{fig:BEVRDAONC}
    \end{subfigure}
    	\par\bigskip
	
    	\begin{subfigure}[t]{0.05\textwidth}
	       \flushright
     \raisebox{.5\height}{\rotatebox{90}{\textbf{RDAO-C}}}
    \end{subfigure}% 
        \begin{subfigure}[t]{0.4\textwidth}
        \centering
        \includegraphics[width=.75\textwidth, height=4cm]{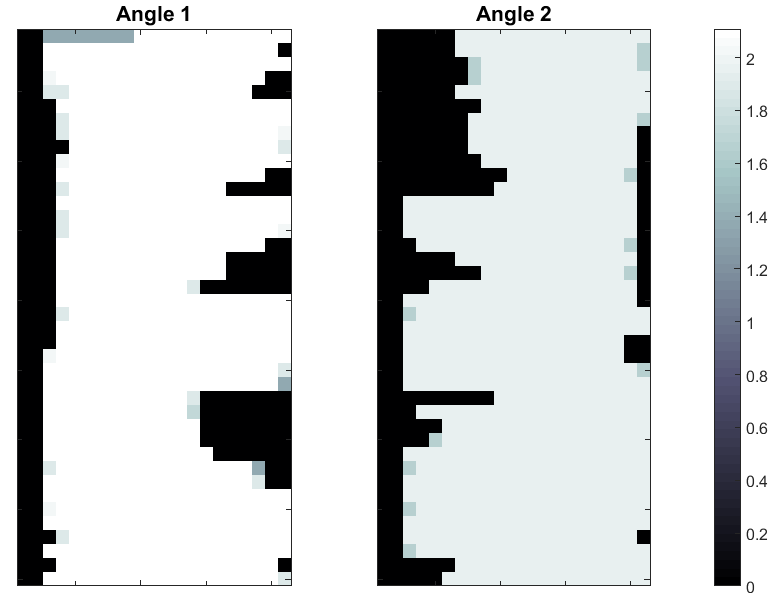}
%        \caption{}%Combined FM of NC DAO}
         \label{fig:FMRDAO}
    \end{subfigure} %	
%         \begin{subfigure}[t]{0.05\textwidth}
%          \hspace{0.05\textwidth}
%     \end{subfigure}% 	 	
       \begin{subfigure}[t]{0.55\textwidth}
        \centering
        \includegraphics[width=.85\textwidth, height=4cm]{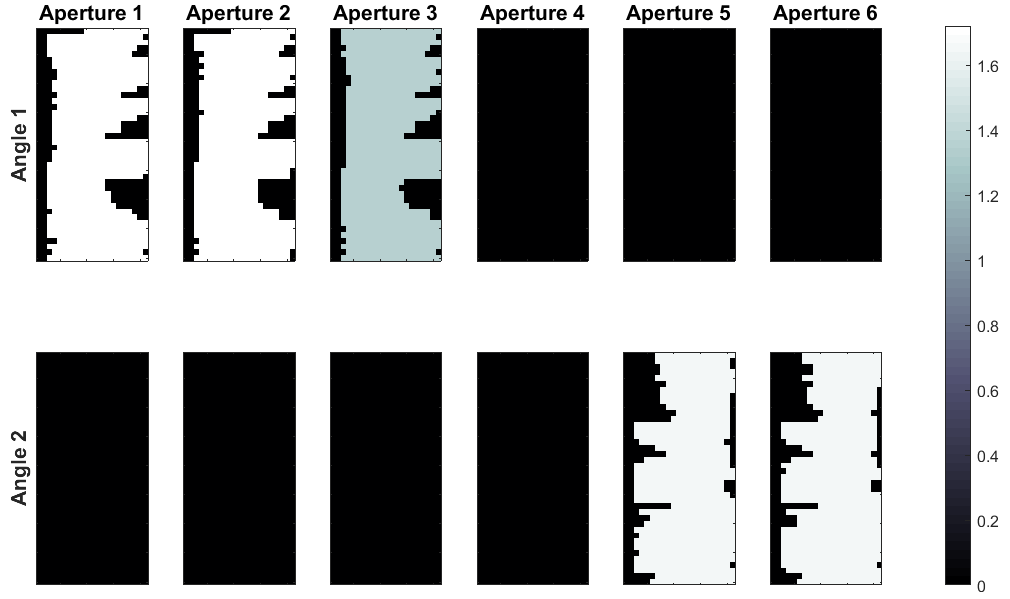}
%        \caption{}%FM of the DAO (no continuity constraints)}
         \label{fig:BEVRDAO}
    \end{subfigure}
    	
    \caption{Aggregated and complete Patient B fluence maps, depicting the best output plans for the four model variants run with a CPG warm start. A log$_{10}$ scale is used in all diagrams to enhance visibility.}
    \label{fig:beamDAO}
    %\vspace{-5pt}
\end{figure}

To better understand the plans, it helps to look at the fluence maps they generate. 
An instance of Patient B's warm-started \eqref{RDAOMod} plans are depicted in Figure~\ref{fig:beamDAO}.
The first column in Figure~\ref{fig:beamDAO} shows the consolidated map of each model variant's plan. This is an aggregated set of fluences, as in Step~1 in Figure~\ref{fig:WarmStartDemo}, however unlike the output of an FMO problem, the figures in the first column are known to decompose into the six fixed, deliverable apertures, depicted in the second column, i.e., a plan. Each row in the figure depicts the output of a different model variant, with the non-robust variants in the first two rows, followed by the robust variants. It is noteworthy that the best incumbent RDAO-C variant uses only 5 of the 6 allowable apertures
This information could be beneficial for clinicians, seeing as fewer apertures reduces setup time. It is unclear, however, whether this reduction would hold at optimality, since none of the models had converged within the three-day time limit.

Model variants with continuity constraints, i.e., DAO-C and RDAO-C, may pose more of a challenge for the solver, but Figure~\ref{fig:beamDAO} provides some insight as to why they are valuable. Consider the non-continuous DAO variant, particularly apertures 4 and 5. Both display undesirable behaviors from Sections~\ref{sec:VertCont} and \ref{sec:HorCont}, with single leaves blocking entire rows, as well as disconnected rows, which would be unlikely to make it through a quality assessment in clinic.
With these behaviors disallowed in DAO-C, the leaf placements are far more desirable.

\begin{table}[]
\centering
\small
    %\color{lightBlue}
\begin{tabular}{@{}clc>{\columncolor[gray]{0.9}}ccccc>{\columncolor[gray]{0.9}}ccccc@{}}
\toprule
 &  &  & \multicolumn{5}{c}{\textbf{Nominal (DAO) Model}} & \multicolumn{5}{c}{\textbf{Robust (RDAO) Model}} \\ \cmidrule(r){4-8} \cmidrule(r){9-13} 
 & \textbf{p} & \textbf{-C} & \multicolumn{1}{c}{$\mathcal{T}^{min}$} & $\mathcal{T}^{ave}$ & $\mathcal{T}^{max}$ & $\mathcal{H}^{ave}$ & $\mathcal{H}^{max}$ & \multicolumn{1}{c}{$\mathcal{T}^{min}$}  & $\mathcal{T}^{ave}$ & $\mathcal{T}^{max}$ & $\mathcal{H}^{ave}$ & $\mathcal{H}^{max}$ \\ \midrule
\multirow{4}{*}{\textbf{Patient A}} & $\textbf{p}_{nom}$ & \xmark & 42.4 & 64.9 & 83.6 & 5.1 & 66.6 & 43.1 & 65.9 & 88.0 & 5.9 & 68.5 \\
 & $\textbf{p}_{real}$ & \xmark & \textbf{40.6} & 64.9 & 83.6 & 5.2 & 66.7 & 42.6 & 65.9 & 88.0 & 6.1 & 68.6 \\
 & $\textbf{p}_{nom}$ & \checkmark & 42.4 & 64.9 & 84.0 & 5.2 & 66.8 & 43.1 & 65.9 & 88.0 & 5.9 & 68.5 \\
 & $\textbf{p}_{real}$ & \checkmark & \textbf{40.7} & 64.9 & 84.1 & 5.4 & 66.9 & 42.6 & 65.9 & 88.0 & 6.1 & 68.6 \\ \midrule
\multirow{4}{*}{\textbf{Patient B}} & $\textbf{p}_{nom}$ & \xmark & 42.4 & 47.0 & 54.3 & 0.5 & 46.0 & 42.4 & 50.2 & 55.6 & 0.6 & 46.2 \\
 & $\textbf{p}_{real}$ & \xmark & \textbf{41.1} & 47.0 & 54.3 & 0.5 & 46.1 & 42.4 & 50.1 & 55.7 & 0.6 & 46.2 \\
 & $\textbf{p}_{nom}$ & \checkmark & 42.4 & 49.3 & 54.5 & 0.6 & 45.7 & 42.6 & 50.5 & 55.6 & 0.6 & 46.8 \\
 & $\textbf{p}_{real}$ & \checkmark & \textbf{41.3} & 49.3 & 54.6 & 0.6 & 45.7 & 42.4 & 50.5 & 55.7 & 0.6 & 46.8 \\ \midrule
\multirow{4}{*}{\textbf{Patient C}} & $\textbf{p}_{nom}$ & \xmark & 42.4 & 63.4 & 79.5 & 2.5 & 63.2 & 42.7 & 63.9 & 81.5 & 2.9 & 67.8 \\
 & $\textbf{p}_{real}$ & \xmark & \textbf{41.4} & 63.3 & 79.5 & 2.4 & 63.5 & 42.4 & 63.8 & 81.5 & 2.8 & 67.9 \\
 & $\textbf{p}_{nom}$ & \checkmark & 42.4 & 63.2 & 80.4 & 2.7 & 65.2 & 42.7 & 63.9 & 81.5 & 2.9 & 67.8 \\
 & $\textbf{p}_{real}$ & \checkmark & \textbf{41.3} & 63.1 & 80.4 & 2.6 & 65.4 & 42.4 & 63.8 & 81.5 & 2.8 & 67.9 \\ \midrule
\multirow{4}{*}{\textbf{Patient D}} & $\textbf{p}_{nom}$ & \xmark & 42.4 & 64.6 & 77.0 & 0.5 & 63.3 & 43.4 & 67.8 & 77.8 & 0.6 & 60.8 \\
 & $\textbf{p}_{real}$ & \xmark & \textbf{41.0} & 64.5 & 77.1 & 0.6 & 66.9 & 43.7 & 67.8 & 77.8 & 0.7 & 63.6 \\
 & $\textbf{p}_{nom}$ & \checkmark & 42.4 & 65.5 & 76.0 & 0.5 & 59.7 & 43.4 & 67.8 & 77.8 & 0.6 & 60.8 \\
 & $\textbf{p}_{real}$ & \checkmark & \textbf{41.3} & 65.4 & 76.0 & 0.6 & 62.9 & 43.7 & 67.8 & 77.8 & 0.7 & 63.6 \\ \midrule
\multirow{4}{*}{\textbf{Patient E}} & $\textbf{p}_{nom}$ & \xmark & 42.4 & 61.9 & 79.4 & 1.1 & 67.3 & 42.4 & 66.2 & 87.3 & 1.6 & 72.9 \\
 & $\textbf{p}_{real}$ & \xmark & \textbf{41.3} & 61.8 & 79.3 & 1.0 & 67.2 & 42.4 & 66.1 & 87.2 & 1.5 & 73.0 \\
 & $\textbf{p}_{nom}$ & \checkmark & 42.4 & 64.3 & 80.1 & 1.2 & 66.4 & 42.4 & 66.2 & 87.3 & 1.6 & 72.9 \\
 & $\textbf{p}_{real}$ & \checkmark & \textbf{40.7} & 64.2 & 80.0 & 1.1 & 66.5 & 42.4 & 66.1 & 87.2 & 1.5 & 73.0 \\ \bottomrule
\end{tabular}
\caption{Unadjusted Warm-Started plans for all five patients with the nominal and with an alternate realized breathing pattern. Underdoses, i.e., $\mathcal{T}_{min} <$ the prescription dose of 42.4~Gy are bolded. }
\label{tbl:robustResults}
\end{table}

In theory, in addition to continuity constraints, clinicians are likely to prefer the delivery guarantees provided by the robust models, making RDAO-C the gold standard for modeling.  
In practice, however, there were some drawbacks when it came to generating robust plans at such large optimality gaps. By nature, robust plans tend to trade off a slight overdose for greater reliability. Since all plans terminated with fairly large optimality gaps, they each had significant overdose, and therefore, a tendency to be fairly robust to uncertainty, impacting the costs versus benefits of using a robust over a nominal model. This phenomenon can be seen illustrated in Table~\ref{tbl:robustResults}, where the unnormalized warm-started \eqref{RDAOMod} plans of the robust and nominal variety are exposed to both nominal and non-nominal breathing patterns, within the uncertainty set $\mathcal{P}$. 
The realized non-nominal breathing pattern $\textbf{p}_{real}$ selected was intended to places more emphasis on the exhale phase for all patients: $\textbf{p}_{real} = \lbrack 0.025, 0.025, 0.125, 0.225, 0.6\rbrack$.

The highlighted columns in Table~\ref{tbl:robustResults} show the minimum dose delivered to a target voxel in a given model variant, given the realized dataset. When the nominal \textbf{p} is used, the target dose of 42.4~Gy is always met, regardless of model. However, when a different \textbf{p} is realized, i.e., $\textbf{p}_{real}$, the hard constraint is no longer met in any of the nominal model variants.

\begin{figure} [htb]
    \centering
    \begin{subfigure}[b]{0.49\textwidth}
        \centering
        \includegraphics[width=\textwidth]{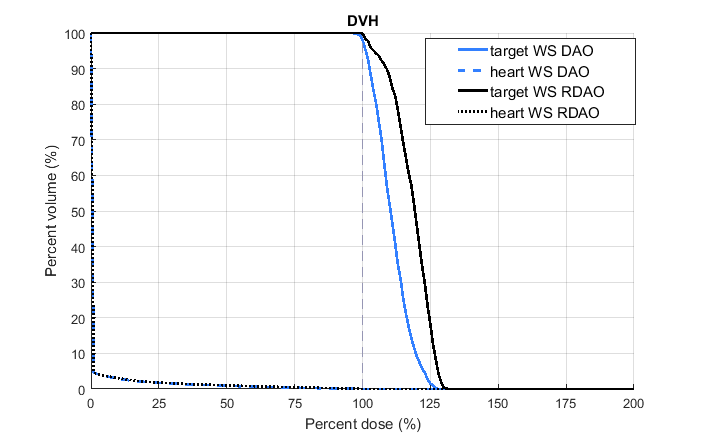}
        \caption{DVH of DAO vs. RDAO WS plan}
        \label{fig:RDAOvDAODVH}
    \end{subfigure}%
    ~~~
    \begin{subfigure}[b]{0.49\textwidth}
        \centering
        \includegraphics[width=\textwidth]{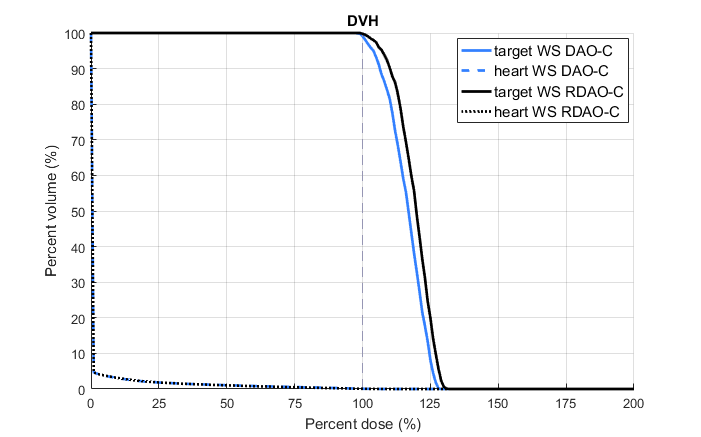}
        \caption{DVH of DAO-C vs. RDAO-C WS plan}
        \label{fig:RDAOCvDAOCDVH}
    \end{subfigure}
    \begin{subfigure}[b]{0.49\textwidth}
        \centering
        \includegraphics[width=.9\textwidth]{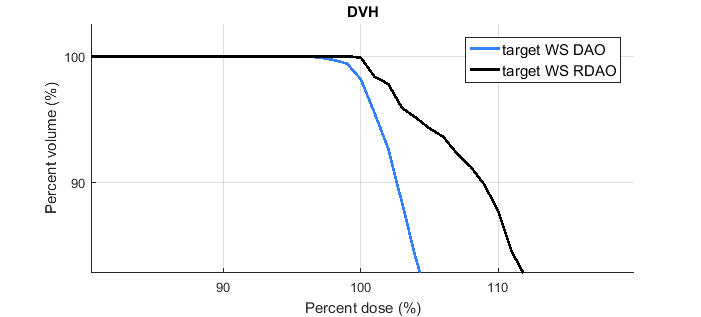}
        \caption{Zoomed-in DAO vs. RDAO WS}
        \label{fig:RDAOvDAOZoomed}
    \end{subfigure}%
    ~~~
    \begin{subfigure}[b]{0.49\textwidth}
        \centering
        \includegraphics[width=.9\textwidth]{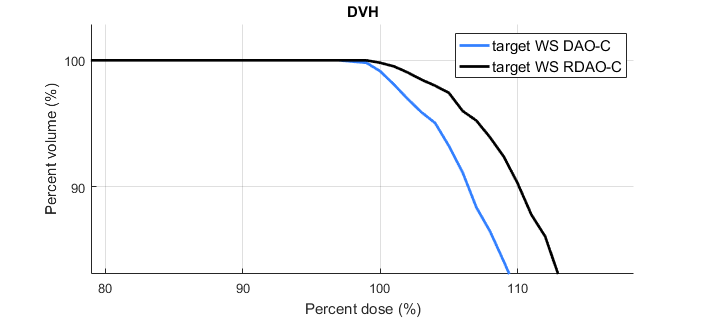}
        \caption{Zoomed-in DAO-C vs. RDAO-C WS}
        \label{fig:RDAOCvDAOCZoomed}
    \end{subfigure}
    \caption{The DVH of Patient B's DAO plans compared to the Robust plans, all generated using the warm start algorithm followed by CPLEX with a $\textbf{p} = \textbf{p}_{real}$ value (rather than nominal). Images \ref{fig:RDAOvDAOZoomed} and \ref{fig:RDAOCvDAOCZoomed} are the zoomed in images of \ref{fig:RDAOvDAODVH} and \ref{fig:RDAOCvDAOCDVH}, respectively. In both images, the lighter blue line is the DAO/-C model, which under-doses in non-nominal $\textbf{p}$ realizations, whereas the black RDAO/-C line continues to deliver 100\% of the dose to 100\% of the target.}
     \label{fig:Pat1RDAOvDAODVHs}
\end{figure}

The drawbacks of the robust method however, can be seen in the remaining columns, where its remaining metrics are consistently outperformed by the DAO and DAO-C models. These observations are reinforced by Figure~\ref{fig:Pat1RDAOvDAODVHs}, which shows the effects of a non-nominal breathing-phase realization on robust and non-robust plans. In figures \ref{fig:RDAOvDAODVH} and \ref{fig:RDAOCvDAOCDVH}, the DAO and DAO-C plans, are visibly more conservative than their robust counterparts. Images \ref{fig:RDAOvDAOZoomed} and \ref{fig:RDAOCvDAOCZoomed} show the cost of these better plans, however, which is the violation of the hard dose constraints.

To sum up these findings, robust and continuous constraints both contain features clinicians want to see, but make the already-difficult DAO MIP model even harder to solve. With the help of the CPG heuristic as a warm start, running a tractable RDAO-C is a more realistic goal than it was previously, in overcoming the incumbent-finding issue. There do, however, remain some hurdles before robust integration lives up to its full potential.

%%%%%%%%%%%%%%%~~~~~  SECTION 5 ~~~~~%%%%%%%%%%%%%%%

\section{Conclusions and Future Work}
\label{Sec:Conclusion} 
This paper proposes a novel, mixed-integer, robust direct aperture optimization (RDAO) model for IMRT treatment planning.
It also provides an original candidate plan generation (CPG) heuristic, which outputs deliverable plans that satisfy the feasibility criteria of the RDAO model.
The planning mechanisms introduced in this paper are demonstrated using five clinical breast cancer case studies with breathing motion uncertainty. 
Three combinations of plan generation techniques are demonstrated: 1) the CPG heuristic alone, 2) the RDAO model alone, and finally, 3) the CPG heuristic as a warm start for the RDAO model. The CPG-RDAO warm start combination was found to dominate both other methods, as expected, as the CPG heuristic was able to provide a boost for the solver, improving its ability to prune sub-optimal solutions. More surprisingly, the CPG heuristic dominated the RDAO model on its own, even when the RDAO model was allowed to run for three days.

While the CPG method did lead to fairly reasonable and deliverable plans, both in standalone and warm-start contexts, the plans did not have the surgical precision that is desired from an optimization methodology, due to the large optimality gaps in the reported solutions. 
This sub-optimality served to undermine the effects of some salient features of the model, such as robustness, which tends to shine in high-precision scenarios. 

To address the solvability concerns, decomposition approaches could be developed to increase the efficiency of the solution procedures for the RDAO model. Column generation could be employed for aperture generation, as in \cite{romeijn2005column}, \cite{men2007exact} and \cite{mahnam2017simultaneous}, or Benders' decomposition could be used to create faster solution algorithms, as in \cite{tacskin2010optimal} and \cite{dursun2019determination}. These improvements would be critical, for handling larger and more complex treatment regions with multiple at-risk organs, as well as the integration of more sophisticated, and potentially non-linear constraints. Additionally, the integration of more clinically relevant robust dose-control metrics (e.g., the CVaR constraints, as in \cite{chan2014robust}) could be added to the base model to improve plan quality, and allow for more detailed exploration of salient trade-offs between robustness and plan quality. 

Another limitation of this work concerns the CPG methodology, which exceeded expectations in terms of generating stand-alone plans, but 
was not powerful enough as a warm start to bring the more precise RDAO model to optimality, which was its intended purpose. 
Going forward, the CPG heuristic may generalize well to other applications, such as VMAT, when it is modeled as a sequence of DAO-style problems at a large number of control-points, as in \cite{men2010ultrafast} and \cite{dursun2019determination}. Alternatively, the CPG heuristic could be used in column generation planning tools, which often benefit from a strong warm starting point.

The models introduced in this paper, and in particular the CPG heuristic, could readily be extended to other application areas. Mathematically, the defining features of the problem at hand are (\textit{i}) a set-covering structure with side constraints, and (\textit{ii}) a coupled subset of binary and continuous variables.
Problems wherein continuous, uniform geometric layers must be preordered and assembled to cover a specified area, e.g., in industries like construction and manufacturing, would mirror the radiation therapy planning process, mathematically. Similarly, difficult 
scheduling problems could be reframed as time-segment selection to cover workloads in applications like worker, or flight shift scheduling, making them a good potential candidate for the CPG heuristic. Finally, the models could be applied to supply chain management, where the robustness could help address demand, or raw-material-arrival uncertainty, and the interdependent continuous-binary relationship would be a good fit for describing facility selection and production-level decisions.

Immediate future directions for this work include targeting structurally similar treatment sites like the liver or lungs, examining methods to speed up RDAO runtime, as well as incorporating new forms of uncertainty and deliverability constraints. For the CPG heuristic, next steps could address some of the limitations of the algorithm, such as its lack of angle flexibility, which led to hard-coded even-angle distribution across apertures in this paper. It would also be of interest to evaluate the CPG's performance against some of the state-of-the art leaf sequencing algorithms both in clinic and from the literature.

\section*{Acknowledgments}
The authors would like to thank the Natural Sciences and Engineering Research Council of Canada for supporting this work. The research was also made possible by the facilities of the Shared Hierarchical Academic Research Computing Network (SHARCNET:www.sharcnet.ca) and Compute/Calcul Canada.

\appendix
 
 \section{Robust Counterpart of \eqref{RobMod} for Polyhedral Uncertainty}
 \label{Sec:robMod}
 Polyhedral uncertainty set \eqref{PDef} leads to 
 infinite possible realizations of $\tilde{\textbf{p}}$ when ${p_i - \underline{p}_i < p_i + \bar{p}_i}$ for at least one phase $i$, and therefore, constraints \eqref{RobModInf} are intractable. \cite{chan2006robust} show that the equivalent robust counterpart of these constraints is both tractable and linear, at the expense of introducing a set of $|\mathcal{V}_T| \times (|\mathcal{I}|+1)$ new dual variables ($y_{i,v}$) and $|\mathcal{V}_T| \times |\mathcal{I}|$ new constraints as follows.
  \begin{subequations}
  \begin{align}
 \nonumber  % \tag{\textbf{RFMO-RC}} \label{RobModRC}\\
  \min & \sum_{s\in \{\cT,\cH\}}\frac{c_s}{|\mathcal{V}_s|}\sum_{v \in \mathcal{V}_s} \sum_{b \in \mathcal{B}} \sum_{i \in \mathcal{I}} p_i D_{v,b,i} \omega_{b} \nonumber\\% \label{NomObj} \\
  \text{s.t.}&\sum_{i \in \mathcal{I}} \Big[ \underline{p}_iy_{0,v} -(\underline{p}_i + \bar{p}_i)y_{i,v} + (p_i - \underline{p}_i)\sum_{b \in \mathcal{B}} D_{b,v,i} \omega_{b} \Big] \geq L_{v} & \forall v \in \mathcal{V}_T, \label{Rob1} \\
  &\sum_{b \in \mathcal{B}}  D_{b,v,i} \omega_{b} - y_{0,v}+ y_{i,v} \geq 0 & \forall i \in \mathcal{I}, v \in \mathcal{V}_T, \label{Rob2} \\
  & \omega_{b} \geq 0  & \forall b \in \mathcal{B}, \label{Rob5}\\
  &y_{i,v} \geq 0 & \forall i \in \mathcal{I}, v \in \mathcal{V}_T, \label{Rob3}\\
  &y_{0,v} ~ \textrm{URS} & \forall v \in \mathcal{V}_T. \label{Rob4}
  \end{align}
  \end{subequations}
 The above robust counterpart to the \eqref{RobMod} model finds the worst-case realization of $\tilde{\textbf{p}}$ for each voxel and optimizes over this realization, thereby immunizing the problem against the worst-case \textit{feasible} realization of the uncertainty set for any chosen intensities.

 \section {Complete Robust Direct Aperture Model}
  \label{Sec:fullMod}
 \definecolor{Gray}{gray}{0.9}
 The following is the complete \eqref{RDAOMod} formulation, with all relevant Section~\ref{Sec:Direct} deliverability constraints.
 The shaded lines are the base model, which assuming $|\mathcal{P}| = 1$ (i.e., no uncertainty) is the FMO model. When the shaded region is run and there is uncertainty ($|\mathcal{P}| > 1$), it is denoted the RFMO model. The remaining unannotated constraints are the classical deliverability constraints. When they are included in the model and certain, the DAO model is being run. If uncertain, the RDAO model is being run. Finally, when the boxed vertical and horizontal continuity constraints are included with the classical deliverability constraints on top of a certain base model, the model is denoted DAO-C, and when uncertain, the model is RDAO-C. 
 \begin{table}[h!]
 \centering 
 \small
 \renewcommand{\arraystretch}{1.5}
 \begin{tabular}{c r l r}
  & \cellcolor{Gray} min & \cellcolor{Gray} $\begin{aligned} \sum_{s\in \{\cT,\cH\}}\frac{c_s}{|\mathcal{V}_s|}\sum_{v \in \mathcal{V}_s} \sum_{b \in \mathcal{B}} \sum_{i \in \mathcal{I}} \sum_{a \in \mathcal{A}} p_i D_{v,b,i} w_{b,a} \end{aligned} $ & \cellcolor{Gray} \\ 
 & \cellcolor{Gray} s.t. & \cellcolor{Gray} $\begin{aligned} \sum_{b \in \mathcal{B}} \sum_{i \in \mathcal{I}} \tilde{p}_i D_{b,v,i} \omega_{b} \geq L_{v} \end{aligned} $ & \cellcolor{Gray} $\forall v \in \mathcal{V}_T, \tilde{\textbf{p}} \in \mathcal{P},$ \\ 
 & &  $\begin{aligned} w_{b,a}  \leq Mx_{b,a} \end{aligned} $ & $\forall b \in \mathcal{B}, a \in \mathcal{A},$ \\ 
 & &  $\begin{aligned} w_{b,a}  \leq f_a + M(1 - x_{b,a}) \end{aligned} $ & $\forall b \in \mathcal{B}, a \in \mathcal{A}, $ \\ 
 & &  $\begin{aligned} w_{b,a}  \geq f_a -M(1 - x_{b,a}) \end{aligned} $ & $  \forall b \in \mathcal{B}, a \in \mathcal{A},$ \\ 
 & &  $\begin{aligned} \sum_{b \in \mathcal{B}_\theta} x_{b,a} \leq |\mathcal{B}_\theta|u_{a,\theta} \end{aligned} $ & $\forall a \in \mathcal{A}, \theta \in \Theta, $ \\ 
 & &  $\begin{aligned} \sum_{\theta \in \Theta} u_{a,\theta} = 1 \end{aligned} $ & $\forall a \in \mathcal{A}, $ \\ 
 & &  $\begin{aligned} l_{q,k+1,\theta,a} \geq l_{q,k,\theta,a} \end{aligned} $ & $  \forall k \in \mathcal{K}', q \in \mathcal{Q}, \theta \in \Theta, a \in \mathcal{A}, $ \\ 
 & &  $\begin{aligned} r_{q,k,\theta,a} \geq r_{q,k+1,\theta,a} \end{aligned} $ & $\forall k \in \mathcal{K}', q \in \mathcal{Q}, \theta \in \Theta, a \in \mathcal{A}, $ \\ 
 & &  $\begin{aligned} x_{q,k,\theta,a} = -1+ l_{q,k,\theta,a}+r_{q,k,\theta,a} \end{aligned} $ & $\forall k \in \mathcal{K}, q \in \mathcal{Q}, \theta \in \Theta, a \in \mathcal{A}, $ \\ 
 %\cline{2-4}%
 \hline
  \multicolumn{1}{|c|}{\multirow{8}{*}{\rotatebox{90}{\textbf{Continuity Constraints (``-C'')~~~}}}} 
 && $\begin{aligned}  j_{q,\theta,a} =  -1 + \bar{j}_{q,\theta,a} + \underline{j}_{q,\theta,a} \end{aligned} $ & \multicolumn{1}{r|}{$\forall q \in \mathcal{Q}, \theta \in \Theta, a \in \mathcal{A}, $} \\ 
  \multicolumn{1}{|l|}{} &&  $\begin{aligned}  j_{q,\theta,a} \leq \sum_{k \in \mathcal{K}} x_{q,k,\theta,a} \end{aligned} $ & \multicolumn{1}{r|}{$ \forall q \in \mathcal{Q}, \theta \in \Theta, a \in \mathcal{A}, $} \\ 
  \multicolumn{1}{|l|}{} &&  $\begin{aligned} |\mathcal{K}| \times  j_{q,a,s} \geq \sum_{k \in \mathcal{K}} x_{q,k,\theta,a} \end{aligned} $ & \multicolumn{1}{r|}{$\forall q \in \mathcal{Q}, \theta \in \Theta, a \in \mathcal{A}, $} \\ 
  \multicolumn{1}{|l|}{} &&  $\begin{aligned} \bar{j}_{q,\theta,a} \leq \bar{j}_{q+1,\theta,a}  \end{aligned} $ & \multicolumn{1}{r|}{$\forall q \in \mathcal{Q}', a \in \mathcal{A}, s \in S, $} \\ 
  \multicolumn{1}{|l|}{} &&  $\begin{aligned} \underline{j}_{q+1,\theta,a} \leq \underline{j}_{q,\theta,a} \end{aligned} $ & \multicolumn{1}{r|}{$\forall q \in \mathcal{Q}', a \in \mathcal{A}, s \in S, $} \\ 
  \multicolumn{1}{|l|}{} &&  $\begin{aligned} j_{q,\theta,a} + j_{q-1,\theta,a} - \sum_{\delta=k+1}^{|\mathcal{K}|} x_{q,\delta,\theta,a} \leq 1 + \sum_{\delta=1}^{k} x_{q-1,\delta,\theta,a}  \end{aligned} $ & \multicolumn{1}{r|}{$\forall k \in \mathcal{K}, q \in \mathcal{Q}'', \theta \in \Theta, a \in \mathcal{A},$}  \\ 
  \multicolumn{1}{|l|}{} & & $\begin{aligned} j_{q,\theta,a} + j_{q-1,\theta,a} - \sum_{\delta=1}^{|\mathcal{K}|-k} x_{q,\delta,\theta,a} \leq 1 + \sum_{\delta=|\mathcal{K}|-k+1}^{|\mathcal{K}|} x_{q-1,\delta,\theta, a} \end{aligned} $ & \multicolumn{1}{r|}{$\forall k \in \mathcal{K}, q \in \mathcal{Q}'', a \in \mathcal{A}, s \in S, $} \\ 
  \multicolumn{1}{|l|}{} &&  $\begin{aligned}  j_{q,\theta,a}, \bar{j}_{q,\theta,a},\underline{j}_{q,\theta,a} \in \{0,1\} \end{aligned} $ & \multicolumn{1}{r|}{$ \forall q \in \mathcal{Q}, \theta \in \Theta, a \in \mathcal{A}, $} \\ 
 %\cline{2-4}%
 \hline
 &\cellcolor{Gray} & \cellcolor{Gray} $\begin{aligned} w_{b,a} \geq 0 \end{aligned} $ & \cellcolor{Gray} $\forall b \in \mathcal{B}, a \in \mathcal{A}, $ \\ 
 & &  $\begin{aligned} f_a  \geq 0 \end{aligned} $ & $\forall a \in \mathcal{A}, $ \\ 
 & &  $\begin{aligned} x_{b,a} \in \{0,1\} \end{aligned} $ & $\forall b \in \mathcal{B}, a \in \mathcal{A}, $ \\ 
 & &  $\begin{aligned} u_{a,\theta}\in \{0,1\}  \end{aligned} $ & $\forall a \in \mathcal{A}, \theta \in \Theta, $ \\ 
 & &  $\begin{aligned} l_{q,k,\theta,a},r_{q,k,\theta,a} \in \{0,1\}  & \end{aligned} $ & $\forall k \in \mathcal{K}, q \in \mathcal{Q}, \theta \in \Theta, a \in \mathcal{A}.  $ \\
 \end{tabular}
 \end{table}

 \section{CPG-S $\alpha$-Parameter Selection}
 \label{app:alphaTune}
\begin{figure}[htb]
      \centering
      \includegraphics[width=.85\textwidth]{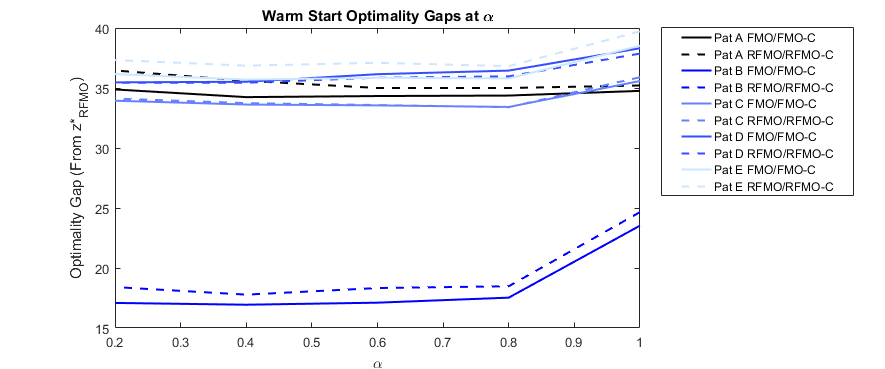}
     \caption{The objective function gap at various $\alpha$ values between 0.2 and 1.}
      \label{fig:WarmStartParamFitting}
  \end{figure}

 The value of the \eqref{WSMod} parameter $\alpha$ was determined by running the problem at different $\alpha$ values in increments of 0.2 and choosing the best plan in terms of objective function value gap. The variation in objective function value as a function of $\alpha$ is depicted in Figure \ref{fig:WarmStartParamFitting}.  The plot starts at $\alpha=0.2$, since the problem collapses back into a fluence map optimization as $\alpha$ approaches 0, generating an optimality gap of 100\% across all cases. Note that the bounds changed very little for $0.2 \leq \alpha \leq 0.8$, so in general, choosing any value in this range should lead to a high-quality warm start.

 \section{A Note on RDAO Uncertainty Levels}
 \label{app:uncertain}
 The RDAO algorithm provides immunization against the worst-case realization of an uncertainty set, which is chosen to be the polyhedral set \eqref{PDef} in this paper, with $\underline{\textbf{p}}=\bar{\textbf{p}}=\textbf{0.1}$. 
 To give some insight into the impact of the magnitude of this uncertainty set, $\underline{\textbf{p}}$ and $\bar{\textbf{p}}$ are varied in Table \ref{tbl:robustPtable}. The objective function values of the CPG-proposed RDAO plan $z^{CPG}$ for each of four different uncertainty conditions, around the nominal $\textbf{p} = \lbrack 0.125, 0.125, 0.125, 0.125, 0.5\rbrack$ are reported. Note that the DAO model (i.e., $\underline{\textbf{p}}=\bar{\textbf{p}}= 0$) and previously studied $\underline{\textbf{p}}=\bar{\textbf{p}} = 0.1$ results are included with the new smaller $\underline{\textbf{p}}=\bar{\textbf{p}} = 0.05$ and larger $\underline{\textbf{p}}=\bar{\textbf{p}} = 0.125$ uncertainty sets for reference.

 \begin{table}[h]
 \centering
 \small
 \begin{tabular}{@{}cccccc@{}}
 \toprule
  & \multicolumn{5}{c}{Patient} \\ \cmidrule(l){2-6} 
 $~\underline{\textbf{p}}=\bar{\textbf{p}~}$ & A & B & C & D & E \\ \midrule
 \textbf{0} & 55.3 & 37.7 & 51.5 & 50.6 & 51.7 \\
 \textbf{0.05} & 56.8 & 38.1 & 52.2 & 51.0 & 52.7 \\
 \textbf{0.1} & 56.8 & 38.2 & 51.8 & 51.3 & 53.0 \\
 \textbf{0.125} & 58.1 & 38.6 & 52.3 & 51.6 & 53.7 \\ \bottomrule
 \end{tabular}
 \caption{The impact robust uncertainty set size on $z^{CPG} $.}
 \label{tbl:robustPtable}
 \end{table}

 The overall trend is that the objective function increases monotonically as the uncertainty set grows larger. This is intuitive, seeing as a plan may have to compensate for larger proportions of time in the worst case phase-realization, thereby giving more dose to the target and healthy organs in order to ensure full dose to target. This trend is not true 100\% of the time in the CPG case (as in Patient C, 0.05 exceeds 0.1) since the CPG heuristic finds non-optimal solutions, however it would be the case at optimality.

\bibliographystyle{unsrtnat}  
\bibliography{references}

\end{document}